\documentclass[ijoc,nonblindrev]{informs3} 

\OneAndAHalfSpacedXII 

\usepackage{amsmath,amsfonts,latexsym,amssymb,amsbsy,tabularx}
\usepackage{mathtools,bm}
\usepackage{mathrsfs}

\usepackage{algorithm}
\usepackage[noend]{algpseudocode}
\usepackage{graphicx}
\usepackage{todonotes}
\usepackage{url}
\usepackage{hyperref}
\usepackage{booktabs}
\usepackage{multirow}
\usepackage{caption}
\usepackage{subcaption}		
\usepackage[inline]{enumitem}

\usepackage{natbib}
 \bibpunct[, ]{(}{)}{,}{a}{}{,}%
 \def\bibfont{\small}%

 \hypersetup{
	colorlinks,
	linkcolor={red!50!black},
	citecolor={blue!50!black},
	urlcolor={blue!80!black}
}

\TheoremsNumberedThrough     

\EquationsNumberedThrough    

\newcommand{\R}{\mathbb{R}}
\newcommand{\Z}{\mathbb{Z}}
\newcommand{\B}{\{0,1\}}
\newcommand{\varIndex}{I}
\newcommand{\domain}{D}
\newcommand{\feasible}{\mathcal{X}}

\newcommand{\bb}{\bm{b}}
\newcommand{\bc}{\bm{c}}
\newcommand{\bx}{\bm{x}}
\newcommand{\by}{\bm{y}}
\newcommand{\bw}{\bm{w}}

\newcommand{\bS}{\bm{S}}

\newcommand{\blambda}{\bm{\lambda}}

\DeclareMathOperator{\conv}{conv} 
\DeclareMathOperator{\proj}{Proj}
\DeclareMathOperator{\projx}{\proj_{\bx}}

\newcommand{\lagsubprob}{\mathcal{L}}
\newcommand{\astar}{$A^*$}

\newcommand{\valueF}{V}
\newcommand{\transitionF}{\phi}
\newcommand{\costF}{g}
\newcommand{\feasibleDP}{X}
\newcommand{\stateSpace}{\mathcal{S}}

\newcommand{\dd}{\mathcal{D}}

\newcommand{\nodes}{\mathcal{N}}
\newcommand{\arcs}{\mathcal{A}}

\newcommand{\maxwidth}{\mathcal{W}}
\newcommand{\width}{w}
\newcommand{\incoming}{\arcs^\textnormal{in}}
\newcommand{\outgoing}{\arcs^\textnormal{out}}
\newcommand{\source}{s}	
\newcommand{\target}{t}	    

\newcommand{\sol}{\textnormal{Sol}}
\newcommand{\pathvar}{\bm{x}}
\newcommand{\pathdd}{p}
\newcommand{\paths}{\mathcal{P}}
\newcommand{\rootnode}{\mathbf{r}}
\newcommand{\terminalnode}{\mathbf{t}}
\newcommand{\rt}{$\rootnode-\terminalnode$}

\newcommand{\arcValue}{v}
\newcommand{\arcLength}{\ell}

\newcommand{\flow}{\mbox{NF}}

\newcommand{\merge}{\oplus}
\newcommand{\bigmerge}{\bigoplus}

\newcommand{\srelaxed}{\widetilde{\bS}}

\newcommand{\Qmin}{Q_{\min}}
\newcommand{\Qmax}{Q_{\max}}

\newcommand{\algTopDownDD}{\textsf{TopDownDD}}
\newcommand{\algDD}{\textsf{ConstructDD}}

\newcommand{\algDDMerge}{\textsf{MergeDDNodes}}
\newcommand{\algDDDiscard}{\textsf{DiscardDDNodes}}

\newcommand{\algReturn}{\textbf{return}}

\begin{document}



\RUNTITLE{Decision Diagrams for Discrete Optimization: A Survey of Recent Advances}

\TITLE{Decision Diagrams for Discrete Optimization: A Survey of Recent Advances}

\RUNAUTHOR{Castro, Cire, and Beck}

\ARTICLEAUTHORS{%
	\AUTHOR{Margarita P. Castro}
	\AFF{Department of Industrial and Systems Engineering, Pontificia Universidad Cat\'olica de Chile, Santiago, Chile,\\ \EMAIL{margarita.castro@ing.puc.cl } \URL{}}
	\AUTHOR{Andre A. Cire}
	\AFF{Department of Management, University of Toronto Scarborough and Rotman School of Management,\\ Toronto, Ontario M1E 1A4, Canada, \EMAIL{andre.cire@utoronto.ca} \URL{}}
	\AUTHOR{J. Christopher Beck}
	\AFF{Department of Mechanical and Industrial Engineering, University of Toronto,
		Toronto, Ontario M5S 3G8, Canada,\\ \EMAIL{jcb@mie.utoronto.ca} \URL{}}
} 

\ABSTRACT{%
	In the last decade, decision diagrams (DDs) have been the basis for a large array of novel approaches for modeling and solving optimization problems. Many techniques now use DDs as a key tool to achieve state-of-the-art performance within other optimization paradigms, such as integer programming and constraint programming. This paper provides a survey of the use of DDs in discrete optimization, particularly focusing on recent developments. We classify these works into two groups based on the type of diagram (i.e., exact or approximate) and present a thorough description of their use. We discuss the main advantages of DDs, point out major challenges, and provide directions for future work.
}



\maketitle

%


\section{Introduction}

Decision diagrams (DDs) are graph-based structures with a large number of applications in computer science and operations research literature. While they have a long history in the Boolean function community \citep{akers1978binary,bryant1986graph,wegener2000branching}, their use in optimization is much more recent. For example, in the past ten years, researchers have successfully applied DDs to methodologies in scheduling \citep{cire2013multivalued}, routing \citep{kinable2017hybrid,castro2019mPDTSP}, and healthcare applications \citep{guo2021logic,riascos2020branch}, to name a few, and related literature and applications are growing steadily.

In the context of optimization, a DD is a directed acyclic graph that encodes solutions and/or their associated costs as directed paths from a root node to a terminal node. The benefit of this representation is that it provides an explicit and potentially compact representation of the solution space of a problem that exposes network structure. Such a network can be manipulated directly, e.g., to obtain valid bounding procedures \citep{bergman2014optimization}, or can be  integrated with other optimization techniques. For example, DDs have been combined with integer programming (IP) solvers using a network-flow formulation over the underlying graph \citep{behle2007binary} and with constraint programming (CP) technologies by developing inference procedures \citep{andersen2007constraint}.

This paper provides an in-depth treatment of DDs for discrete optimization with a focus on recent advances. Our objective is to provide a systematic classification of the field, identify unifying themes, and discuss challenges and opportunities that may be the basis of new research. Specifically, our classification partitions the area into the two pervasive DD types in related works, i.e., exact and approximate DDs. Exact DDs encode the exact problem, in that any property that holds for the DD is also valid for the original problem. Approximate DDs, in contrast, provide an over- or under-approximation of the feasible space or the objective function, and are the basis of combined and enumerative approaches  (e.g., \citealt{bergman2016discrete}).


Using this classification, we divide works employing exact DDs into four themes based on the methodology's purpose: (i) modeling,  (ii) feasibility checking, (iii) solution extraction, and (iv) solution-space analysis. Table \ref{tab:lit_exact} presents the paper classification for exact DDs and its respective sub-classes. Similarly, we divide the works that consider approximate DDs into three themes: (i) approximate DD compilation, (ii) DD-based bounds, and (iii) CP propagation. Table \ref{tab:lit_approx} summarize the works that focus on approximate DDs for each class and sub-class.

\begin{table}[htbp]
	\scriptsize 
	\centering
	\caption{Paper classification for exact decision diagrams.}
	\begin{tabular}{ll}
		\toprule
		Sub-class & Papers \\
		\midrule
		\multicolumn{2}{c}{Modeling} \\
		\midrule
		Recursive Model & \cite{hooker2013decision}, \cite{hooker2017job}, \cite{bergman2016discrete}.  \vspace{0.5em} \\
		\multirow{5}[0]{*}{Network Flow Formulation} &  \cite{behle2007binary}, \cite{bergman2016decomposition}, \cite{latour2017combining}, \\
		&  \cite{haus2017scenario}, \cite{bergman2018nonlinear}, \cite{latour2019stochastic}, \\
		&  \cite{serra2019last}, \cite{hosseininasab2019exact}, \cite{cire2019network}, \\
		&  \cite{ploskas2019heat}, \cite{bergman2019binary}, \cite{lozano2020consistent}, \cite{lozano2020decision}, \\
		& \cite{nadarajah2020network}, \cite{bergman2021quadratic}, \cite{mehrani2021models}. \vspace{0.5em} \\
		\multirow{5}[0]{*}{Global Constraints} & \cite{andersen2007constraint}, \cite{cheng2008maintaining}, \cite{cheng2010mdd}, \\
		&  \cite{perez2014improving}, \cite{amilhastre2014compiling},  \cite{perez2015relations},  \\
		&  \cite{perez2015efficient}, \cite{perez2016constructions}, \cite{roy2016enforcing},  \\
		& \cite{perez2017soft}, \cite{perez2018parallel}, \cite{verhaeghe2018compact},  \\
		& \cite{vion2018mdd}, \cite{verhaeghe2019extending}, \cite{de2019compiling}. \vspace{0.5em} \\
		Continuous Variables & \cite{davarnia2021strong}, \cite{salemi2021structure}. \\
		\midrule
		\multicolumn{2}{c}{Feasibility Checking} \\
		\midrule
		General & \cite{nishino2015bdd}, \cite{xue2019embedding}. \vspace{0.5em} \\
		\multirow{2}[0]{*}{Cutting Planes} & \cite{becker2005bdds}, \cite{behle2007binary}, \cite{tjandraatmadja2019target}, \\
		&  \cite{davarnia2020outer}, \cite{castro2021combinatorial}. \vspace{0.5em} \\
		\multirow{2}[0]{*}{Benders Decomposition} & \cite{lozano2018binary}, \cite{guo2021logic}, \\
		& \cite{salemi2021structure}, \cite{bergman2021quadratic}. \vspace{0.5em} \\
		\multirow{2}[1]{*}{Inference} & \cite{subbarayan2008efficient},  \cite{hadzic2009enhanced}, \cite{gange2011mdd},  \\
		& \cite{gange2013explaining}, \cite{kell2015clause}, \cite{jung2021checking}. \\
		\midrule
		\multicolumn{2}{c}{Solution Extraction} \\
		\midrule
		General & \cite{hadzic2004fast}. \vspace{0.5em} \\
		\multirow{2}[1]{*}{Column Generation} & \cite{morrison2016solving}, \cite{kowalczyk2018branch}, \\
		& \cite{raghunathan2018seamless}, \cite{riascos2020branch}. \\
		\midrule
		\multicolumn{2}{c}{Solution-Space Analysis} \\
		\midrule
		Post-Optimiality Analysis & \cite{hadzic2006postoptimality}, \cite{hadvzic2007cost}, \cite{serra2019compact}. \vspace{0.5em} \\
		\multirow{2}[1]{*}{Solution Enumeration} & \cite{bergman2016multiobjective},  \cite{haus2017compact}, \\
		& \cite{suzuki2018fast}, \cite{suzuki2018exact}, \cite{bergman2021network}. \vspace{0.5em} \\
		Polyhedral Analysis & \cite{behle2007facet}, \cite{tjandraatmadja2019target}. \\
		\bottomrule
	\end{tabular}%
	\label{tab:lit_exact}%
\end{table}%

\begin{table}[tbp]
	\scriptsize
	\centering
	\caption{Paper classification for approximate decision diagrams.}
	\begin{tabular}{ll}
		\toprule
		Sub-class & Papers \\
		\midrule
		\multicolumn{2}{c}{Approximate DD Compilation} \\
		\midrule
		Top-down and  & \cite{hadzic2008approximate}, \cite{bergman2011manipulating}, \cite{bergman2014heuristic},  \\
		Iterative Refinement & \cite{frohner2019merging}, \cite{frohner2019towards}, \\
		& \cite{de2020single}, \cite{horn2021based}. \vspace{0.5em} \\
		Other Construction & \cite{cire2014separation}, \cite{bergman2016theoretical},\cite{bergman2017finding}, \\
		Algorithms & \cite{romer2018local}, \cite{horn2021based}, \cite{horn2021common}. \vspace{0.5em} \\
		\multirow{2}[1]{*}{Variable Ordering} & \cite{behle2008threshold}, \cite{bergman2011manipulating}, \cite{cappart2019improving}, \\
		& \cite{karahalios2021variable}, \cite{parjadis2021improving}. \\
		\midrule
		\multicolumn{2}{c}{DD-based Bound} \\
		\midrule
		\multirow{4}[1]{*}{Dual Bounds} & \cite{andersen2007constraint}, \cite{cire2013multivalued}, \cite{kell2013mdd},  \\
		& \cite{bergman2014optimization}, \cite{hooker2017job}, \cite{kinable2017hybrid},  \\
		& \cite{van2018multi}, \cite{maschler2018multivalued}, \cite{castro2019relaxedbdds},  \\
		& \cite{castro2019mPDTSP}, \cite{castro2020SolvingDF}, \cite{van2020graph}, \cite{van2021graph}.  \vspace{0.5em} \\
		\multirow{2}[0]{*}{Lagrangian Bounds} & \cite{bergman2015improved}, \cite{bergman2015lagrangian}, \cite{hooker2019improved},  \\
		& \cite{castro2019mPDTSP}, \cite{tjandraatmadja2020incorporating}, \cite{lange2021efficient}. \vspace{0.5em} \\
		\multirow{2}[0]{*}{Primal Bounds} & \cite{kell2013mdd}, \cite{bergman2014heuristic}, \\
		& \cite{o2019decision}, \cite{horn2019decision}. \vspace{0.5em} \\
		\multirow{2}[1]{*}{Branch-and-Bound} & \cite{bergman2014parallel}, \cite{bergman2016discrete}, \cite{gonzalez2020bdd},  \\
		& \cite{gonzalez2020integrated}, \cite{gillard2021improving}. \\
		\midrule
		\multicolumn{2}{c}{CP Propagation} \\
		\midrule
		& \cite{andersen2007constraint}, \cite{hadzic2008propagating}, \cite{hoda2010systematic},  \\
		& \cite{hadzic2009enhanced}, \cite{cire2012mdd}, \cite{bergman2014mdd}, \\
		& \cite{perez2017mddsampling},  \cite{perez2017mdds}, \cite{kinable2017hybrid}. \\
		\bottomrule
	\end{tabular}%
	\label{tab:lit_approx}%
\end{table}%

The paper is organized as follows. Section \ref{sec:preliminaries} provides a background on DDs and the notation used throughout the paper. Section \ref{sec:exactDD} presents works related to exact DDs, highlighting recent advances in the field. Similarly, Section \ref{sec:approxDD} focuses on approximate DDs and novel methodologies in the literature. We concluded this survey in Section \ref{sec:conclusions} with final remarks and future work directions. 

\section{Preliminaries} \label{sec:preliminaries}

We now formally define DDs and present the notation used throughout the paper. We start with the DD encoding of optimization problems in Section \ref{sec:representation}, show the connections to recursive models in Section \ref{sec:pre_recursive}, and discuss basic approximation concepts in Section \ref{sec:pre_approx}. While there are several DD variants in the literature (see, e.g., \citealt{wegener2000branching}), we primarily focus on ordered decision diagrams. That is, layered diagrams where each layer is associated with a single variable. This DD variant is the most commonly used by the discrete optimization community and, as such, most contributions are based on this structure.

\subsection{DD Representation of Optimization Problems}
\label{sec:representation}

In this paper, we consider maximization problems \ref{model:dop_general} of the form
\begin{align}
	\max_{\bx} & \big\{ f(\bx) \; : \; \bx\in \feasible \; \big \}. \tag{$\mathscr{P}$} \label{model:dop_general}
\end{align}
where $\bx = (x_1, x_2, \dots, x_n)$ is a tuple of $n$ variables, $\feasible \subset \mathbb{Z}^n$ is a finite (integer) feasible space, 
and $f: \Z^n \rightarrow \R$ is the objective function. Moreover, we let $\varIndex \equiv \{1,\dots,n\}$ denote the variable index set, and consider that each $x_i$ assumes values in a finite domain $\domain_i \subset \mathbb{Z}$ where $\feasible \subseteq \domain_1\times\dots\times \domain_n$. 

\begin{figure}[htb]
	\centering
	\includegraphics{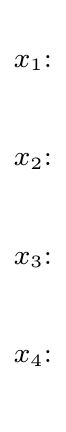}
	\includegraphics{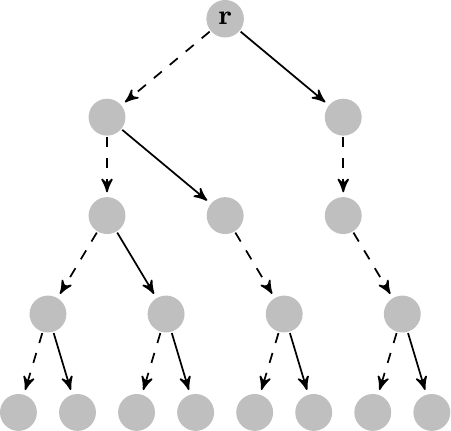}
	\hspace{7em}
	\includegraphics{figs-bdd-vars.pdf}
	\includegraphics{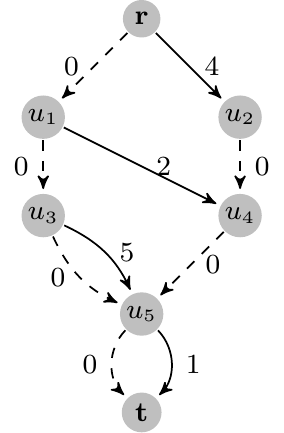}
	\includegraphics{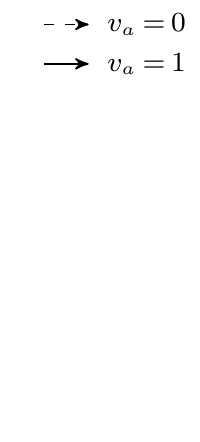}
	
	\caption{A decision tree and a DD for $\feasible=\{ \bx\in \B^4: 7x_1 + 5x_2 + 4x_3 +x_4 \leq 8 \}$.} \label{fig:knap_dd} 
\end{figure}

DDs are graphical representations of solutions to \ref{model:dop_general}, where layers map to variables and arcs map to variable-value assignments (see, e.g., Figure \ref{fig:knap_dd}). Formally, a DD $\dd= (\nodes, \arcs)$ is a layered-directed acyclic graph with node set $\nodes$ and arc set $\arcs$. The node set $\nodes$ is partitioned into $n+1$ layers $\nodes = (\nodes_1,\dots, \nodes_{n+1})$. The first and last layers are the singletons $\nodes_1=\{\rootnode\}$ and $\nodes_{n+1}=\{\terminalnode\}$, respectively, where $\rootnode$ is the root node and $\terminalnode$ is the terminal node. An arc $a = (u,u') \in \arcs$ has a source node $\source(a) = u$ and a target node $\target(a) = u'$ in consecutive layers, i.e., $u' \in \nodes_{i+1}$ whenever $u \in \nodes_i$ for every $i \in \varIndex$.

The solutions in $\feasible$ are mapped to paths in the graph as follows. The arcs emanating from layer $\nodes_i$, $i \in I$, are associated with values in the domain of variable $x_i$. Every arc $a \in \arcs$ with source $\source(a) \in \nodes_i$ has a value $\arcValue_a \in \domain_i$, and each node  $u \in \nodes_i$ has at most $|D_i|$ emanating arcs, each with a different value. Given an arc-specified \rt\ path $\pathdd=(a_1,\dots,a_n)$ with $\source(a_1) = \rootnode$ and $\target(a_n)=\terminalnode$, we let $\pathvar^\pathdd = (\arcValue_{a_1}, \arcValue_{a_2}, \dots, \arcValue_{a_n}) \in \domain_1\times \dots \times \domain_n$ be the $n$-dimensional point encoded by path $\pathdd$. Then, the solutions represented by the DD is
\[\sol(\dd) \equiv \bigcup_{\pathdd \in \paths} \{ \pathvar^\pathdd \},\]

\noindent where $\paths$ is the set of all \rt\ paths in $\dd$. Finally, each arc $a \in \arcs$ is also associated with an arc length $\arcLength_a$. The length encodes, e.g., the contribution to the objective function of paths crossing that arc. We denote the length of a path $\pathdd \in \paths$ by $\arcLength(\pathdd) \equiv \sum_{a\in p} \arcLength_a$. 

We say that $\dd$ exactly represents $\feasible$ if $\sol(\dd) = \feasible$ and $f(\pathvar^\pathdd) = \arcLength(\pathdd)$ for every $\pathdd \in \paths$, i.e., there is a one-to-one correspondence between points in $\feasible$ and paths in $\paths$, and path lengths evaluate to the objective value of their associated solution.

\medskip
\begin{example}\label{exa:knap}
	Consider a knapsack problem with feasible set $\feasible=\{ \bx\in \B^4: 7x_1 + 5x_2 + 4x_3 +x_4 \leq 8 \}$ with weight vector $\bw=(7,5,4,1)$ and a linear objective $f(\bx) = \bc^\top\bx$ with cost vector $\bc=(4,2,5,1)$. Figure \ref{fig:knap_dd} illustrates the set of feasible solutions in $\feasible$ over two graphs. Dashed arrows represent arcs associated with a zero-value variable assignment  (i.e., $\arcValue_a = 0$) and solid arrows correspond to arcs with a one-value variable assignments (i.e., $\arcValue_a = 1$). The lengths of solid and dashed arcs in layer $i$ are $c_i$ and 0, respectively.
	
	The left graph in Figure \ref{fig:knap_dd} corresponds to a search tree (e.g., in a branch-and-bound procedure) where paths from the root $\rootnode$ to each leaf node correspond to feasible variable assignments. The right graph illustrates a DD for $\feasible$. Each DD layer is associated with a variable and arcs denote feasible variable-value assignments. Note that there is a one-to-one correspondence between root-to-leaf paths in the search tree and \rt\ paths in the DD, i.e., the DD exactly represents $\feasible$. Moreover, by construction, each path length evaluates to the objective value of the corresponding solution. Note that any longest path, here defined by $\bx^* = (0, 0, 1, 1)$ with a length of 6, is an optimal solution to the problem.
	\hfill $\blacksquare$
\end{example}

The DD variant above is the most commonly used by the discrete optimization community due to its simpler layer-variable structure. The original DD definition for Boolean functions \citep{bryant1986graph, bryant1992symbolic} differs in that:
\begin{enumerate*}[label=(\alph*)]
	\item nodes are directly associated with variables (i.e., literals) without the layered constraint,
	\item arcs may traverse multiple layers, and
	\item the DD has two terminal nodes, one for feasible and one for infeasible solutions.
\end{enumerate*}
While some researchers employ the classic DD definition to create smaller diagrams for specific applications at the cost of more intricate implementations \citep{perez2014improving}, most works in optimization adopt the restricted DD structure above. A second notable DD variant is a zero-suppressed decision diagram (ZDD), which can compactly represent combinatorial sets and has been used, for instance, in column generation algorithms \citep{morrison2016solving,kowalczyk2018branch} and enumeration procedures \citep{suzuki2018fast}. We refer the reader to the work of \cite{minato1993zero} for further details on ZDDs.

\subsection{Decision Diagrams and Recursive Formulations}
\label{sec:pre_recursive}

A common practice in the literature is to conceptualize DDs based on a recursive formulation (RF) of \ref{model:dop_general}. As investigated by \cite{hooker2013decision}, there is a strong relationship between DDs and the state transition graph in the dynamic programming (DP) literature. DPs represent optimization problems via a recursive model where decisions are made sequentially \citep{bertsekas1995dynamic}. The state transition graph can be obtained by unfolding the recursive model for any possible state in the system, where nodes map to states, arcs represent the state-transition function, and arc lengths encode the immediate rewards.

To detail this relationship, we start by introducing the general form of a recursive model using syntax from the DP literature \citep{bertsekas1995dynamic}. In particular, such models are defined in terms of stages, where transitions across stages are driven by actions (or controls). We consider an $n+1$ stage system with actions $\bx \in \Z^n$, one per stage $i \in \varIndex$. For a given stage $i\in \{1,\dots,n+1\}$, state variables $\bS \in \stateSpace_i$ represent a summary of past actions that reach that state, where the state space $\stateSpace_i$ denotes the set of states reachable at the $i$-th stage. The set of feasible actions $x_i$ at a state $\bS\in \stateSpace_i$ is given by the action set $\feasibleDP_i(\bS)$. 

The system transitions to a new state according to the current state and action that has been applied. Transitions are governed by the transition function $\transitionF_i : \stateSpace_i\times \feasibleDP_i \rightarrow \stateSpace_{i+1}$, which maps the current state-action pair to a state in the next stage. Moreover, each state $\bS  \in \stateSpace_i$ and decision variable $x\in \feasibleDP_i(\bS)$ at stage $i \in \{1,\dots,n+1\}$ incurs an immediate reward $\costF_i(\bS, x)$. 

The optimal actions of a DP model solve the Bellman equation
\begin{align}
	\valueF_i(\bS) = & \max_{x \in \feasibleDP_i(\bS)}\left\{ \costF_i(\bS, x) + \valueF_{i+1}(\transitionF_i(\bS, x) )  \right\}, \qquad \forall i \in \varIndex,     \tag{RF} \label{model:gen_rec} 
\end{align}
where $\valueF_i(\cdot)$ is the value function with $\valueF_{n+1}(\bS) = 0$ for any $\bS \in  \stateSpace_{n+1}$. We assume the initial state is the singleton $\stateSpace_1= \{ \bS_1 \}$, where $\bS_1$ is the root state.

A recursive model is a reformulation of problem \ref{model:dop_general} if
\begin{enumerate}[label=(\alph*)]
	\item There is an one-to-one mapping between solutions $x \in \feasible$ and a state-action trajectory $(\bS_1, x_1), (\bS_2, x_2), \dots, (\bS_n, x_n)$, where $\bS_{i+1} = \transitionF_i(\bS_{i},x_i)$. 
	
	\item For every solution $\bx \in \feasible$ and associated state trajectory 
	$(\bS_1, x_1), (\bS_2, x_2), \dots, (\bS_n, x_n)$, we have $f(\bx) = \sum_{i=1}^n \costF_i(\bS, x_i)$, i.e., the solution value matches the sum of rewards.
\end{enumerate}

\medskip
\begin{example} \label{exa:knap_dp}
	We depict the classical recursive model of the knapsack problem introduced in Example \ref{exa:knap}. The actions are $\bx \in \B^4$ and the state $\bS = Q \subseteq \mathbb{Z}$ captures the load of the knapsack at each stage. For an initial state $Q_1=0$, 
	the Bellman equation is
	\begin{align}
		\valueF_i( Q) = & \max_{x \in \feasibleDP_i(Q)}\left\{ c_i x + \valueF_{i+1}( Q + w_ix )  \right\},  \qquad \forall i \in \{1,\dots,5\}.    \tag{R-KNP} \label{model:knap_rec} 
	\end{align}
	
	The transition function $\transitionF_i(Q,x) = Q + w_ix$ updates the load of the knapsack, while the immediate reward $\costF_i(Q,x) = c_i x$ corresponds to the gain from choosing item $i$. Since the weight of each item is positive, the action space for state $Q \in \stateSpace_i$ is 
	$\feasibleDP_i(Q) = \{ x \in \B: Q + w_ix \leq 8 \}$, for all stages $i \in \{1,\dots,4\}$. \hfill $\blacksquare$
\end{example}

\medskip 
A DD can be obtained immediately from the state-transition graph of \ref{model:gen_rec} \citep{bergman2011manipulating}. Specifically, with each state $\bS \in \stateSpace_i$ we associate a node $u_{\bS} \in \nodes_i$, $i \in \varIndex \cup \{n+1\}$; note that the root node $\rootnode$ corresponds to the initial state $\bS_1$. There exists an arc $a$ with value $\arcValue_a = x$ between nodes $u_{\bS_i} \in \nodes_i$ and $u_{\bS_{i+1}} \in \nodes_{i+1}$, $i \in \varIndex$, if and only if $\bS_{i+1} = \transitionF_i(\bS_{i},x)$. The length of such an arc is the immediate reward of the transition, i.e., $\arcLength_a = \costF_i(\bS_i, x)$. Finally, the nodes in the last layer, in case several terminal states are present, are merged into the single terminal node $\terminalnode$. 

\subsection{DD Approximations} 
\label{sec:pre_approx}

While the size of the DD can be considerably reduced using the classical reduction procedure by 
\cite{bryant1992symbolic} (i.e., by merging isomorphic subgraphs), an exact DD $\dd$ for $\feasible$ is, in general, exponentially large in the number of variables $n$. An alternative to overcome this limitation is to consider a DD approximation instead, i.e., a DD that either under or over-approximates the set of feasible solutions and/or the objective function. 

This idea was first proposed by \cite{andersen2007constraint}, who over-approximates $\feasible$ using a relaxed DD, i.e., a DD $\dd$ with solution set such that $\feasible \subseteq \sol(\dd)$. Thus, every point in $\feasible$ maps to a \rt\ path in $\dd$, but the converse does not necessarily hold. \cite{bergman2018nonlinear}, \cite{hooker2013decision}, and \cite{hooker2019improved} also consider variants that relax the objective function value, in that any \rt\ path $\pathdd$ satisfies $f(\pathvar^\pathdd) \le \arcLength(\pathdd)$. Relaxed DDs provide a discrete relaxation of $\feasible$ that can be used, for instance, to compute optimistic bounds for an optimization problem and have been key in state-of-the-art DD methodologies.

Conversely, \cite{bergman2014heuristic} introduce the concept of restricted DDs, i.e., $\dd$ under-approximates $\feasible$ by capturing only a subset of the feasible solutions. Such an approximation can be obtained, for example, by limiting the number of nodes in each layer and heuristically discarding nodes. The main advantage of this technique is to obtain primal bounds for an optimization problem, as 
we will further discuss in Section \ref{sec:approx_bounds}.

\begin{figure}[htb]
	\centering
	\includegraphics{figs-bdd-vars.pdf}
	\includegraphics{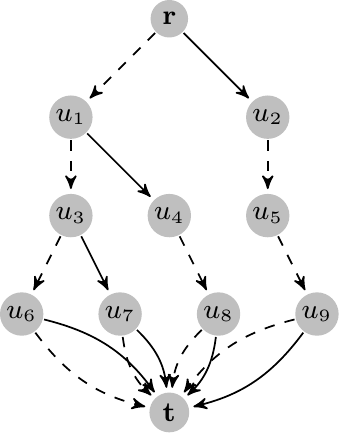}
	\hspace{4em}
	\includegraphics{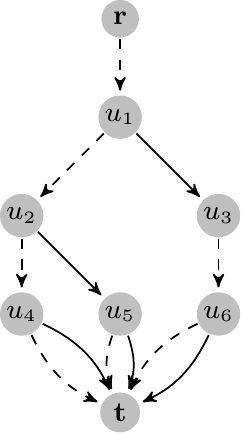}
	\hspace{4em}
	\includegraphics{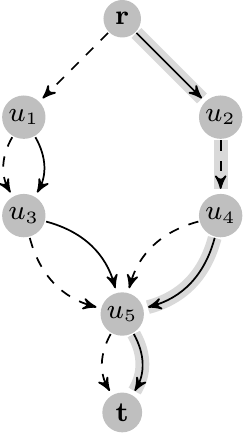}
	\includegraphics{figs-bdd-arc-values.pdf}
	\caption{An exact DD $\dd_1$ (left), a restricted DD $\dd_2$ (middle), and a relaxed DD $\dd_3$ (right) for $\feasible=\{ \bx\in \B^4: 7x_1 + 5x_2 + 4x_3 +x_4 \leq 8 \}$.} \label{fig:knap_dd_all} 
\end{figure}

\medskip
\begin{example} \label{exa:knap_dds_all}
	We recall the knapsack instance from Example \ref{exa:knap} with feasible set $\feasible=\{ \bx\in \B^4: 7x_1 + 5x_2 + 4x_3 +x_4 \leq 8 \}$ and recursive model \ref{model:knap_rec}. Figure  \ref{fig:knap_dd_all} depicts an exact, a relaxed, and a restricted DD for $\feasible$. The paths of the exact DD $\dd_1$ (left diagram) have a one-to-one relationship with the points in $\feasible$. Moreover, each node in $\dd_1$ is associated with a single state of  \ref{model:knap_rec}, e.g., node $u_1$ has $Q=0$ and node $u_2$ has $Q=7$.
	
	The restricted DD $\dd_2$ (middle) represents only a subset of the solutions in $\feasible$. For instance, the solution $(1,0,0,1) \in \feasible$ is not in $\sol(\dd_2)$. This restricted DD can be obtained by discarding the states of \ref{model:knap_rec} that have value greater than $6$. Lastly, relaxed DD $\dd_3$ represents all feasible solutions in $\feasible$ and some infeasible assignments. Specifically, the shaded \rt\ path in $\dd_3$ corresponds to infeasible solution $(1,0,1,1)\notin \feasible$. \hfill $\blacksquare$
\end{example}
\medskip

When generating decision diagrams from a recursive model \ref{model:gen_rec}, each node and arc in an exact DD maps to a state and action, respectively. Thus, the transition and reward functions are well-defined in such settings, as they can be inherited directly from \ref{model:gen_rec}. In a relaxed DDs, however, this  property does not hold because relaxing a DD modifies the state space by incorporating infeasible paths into the solution set. 

\cite{hooker2017job} provides a formal connection between \ref{model:gen_rec} and a relaxed DD through the use of \textit{merged states}. More specifically, given two states $\bS,\bS' \in \stateSpace_i$, we define $\srelaxed=\bS\merge\bS'$ as a merged state where $\merge$ is an appropriate merging operator, i.e.,  $\merge$ satisfies the following properties for any $i \in \varIndex$:
\begin{enumerate}[label=(C\arabic*)]
	\item The set of feasible actions over states $\bS$ and $\bS'$ are also feasible over $\srelaxed$, i.e.,  $\feasibleDP_i(\bS), \feasibleDP_i(\bS') \subseteq \feasibleDP_i(\srelaxed)$. \label{cond:rs_subset}
	\item The immediate reward at state $\srelaxed$ is greater than or equal to the immediate reward at states $\bS$ and $\bS'$. Thus, for any $x \in \feasibleDP_i(\bS)$ and  $x' \in \feasibleDP_i(\bS')$ we have $\costF_i(\bS, x) \leq \costF_i(\srelaxed, x)$ and $\costF_i(\bS', x') \leq \costF_i(\srelaxed, x')$, respectively. \label{cond:rs_cost}
\end{enumerate}

Following the above conditions, we say that $\srelaxed$ relaxes a state $\bS \in \stateSpace_i$ if  $\feasibleDP_i(\bS)\subseteq \feasibleDP_i(\srelaxed)$ and the immediate reward function over any $x\in \feasibleDP_i(\bS)$ is larger for $\srelaxed$, i.e., $\costF_i(\bS, x) \leq \costF_i(\srelaxed, x)$. Properties \ref{cond:rs_subset} and \ref{cond:rs_cost} are necessary (but not sufficient) conditions to define a proper relaxation of \ref{model:gen_rec}. \cite{hooker2017job} shows that operator $\merge$ defines a proper relaxation for \ref{model:gen_rec} if, in addition to \ref{cond:rs_subset} and \ref{cond:rs_cost}, we impose a condition over the transition function:
\begin{enumerate}[label=(C\arabic*)]
	\setcounter{enumi}{2}
	\item If $\srelaxed$ relaxes state $\bS \in \stateSpace_i$, then, given any value $x\in \feasibleDP_i(\bS)$, $\transitionF_i(\srelaxed, x)$ relaxes state $\transitionF_i(\bS, x)$, for all $i \in \varIndex$. Thus, $\srelaxed$ defines a relaxed state in the following stage for each feasible action.  \label{cond:rs_transition}
\end{enumerate}

The relevance of establishing such a connection is to provide a framework for building relaxed DDs using \ref{model:gen_rec}. For example, one could construct a relaxed DD using a top-down approach starting with $\rootnode$ and building one layer at a time. If the number of nodes in a layer (width) is too large, nodes can be merged using a state redefinition that satisfies \ref{cond:rs_subset}-\ref{cond:rs_transition}. We provide details in Section \ref{sec:approx_construction}.

\begin{example} \label{exa:knap_rstate}
	Recall the knapsack instance from Example \ref{exa:knap} with feasible set $\feasible=\{ \bx\in \B^4: 7x_1 + 5x_2 + 4x_3 +x_4 \leq 8 \}$ and recursive model \ref{model:knap_rec}. We define the merge operator over states $Q, Q' \in \stateSpace_i$ at stage $i \in \{1,\dots,4\}$ as the minimum over both quantities, i.e.,	$ Q\merge Q' = \min\{Q, Q'\}$.	Operator $\merge$ satisfies \ref{cond:rs_subset} since any feasible solution for a knapsack load $Q$ or $Q'$ is also feasible for $\min\{Q, Q'\}$. Condition \ref{cond:rs_transition} also holds for $\merge$ since the transition function is an increasing function over $Q$ (see Example \ref{exa:knap_dp}). Lastly, \ref{cond:rs_cost} holds since the immediate reward $\costF_i(Q,x)= c_ix$ is independent of the current state $Q$. Note that $\dd_3$ in Figure \ref{fig:knap_dd_all} is a relaxed DD created with this merging operator by choosing to merge nodes with small state value first. \hfill $\blacksquare$
\end{example}

\section{Exact Decision Diagrams}\label{sec:exactDD}

Exact DDs encode the set of solutions of a discrete optimization problem as paths over a directed acyclic graph. As shown in Table \ref{tab:lit_exact}, we distinguish four different purposes observed in the literature for constructing an exact DD: modeling, feasibility checking, solution extraction,  and solution-space analysis. 

We review each of these purposes and show how to integrate them into integer programming (IP) and constraint programming (CP) solvers. This section focuses on works where one or multiple exact DDs represent either the complete problem or a subset of its constraints, i.e., where no merging operation is applied. We first present different modeling techniques based on DDs.  We then review works that use DDs to check feasibility or to  extract solutions. The last subsection describes enumeration procedures and post-optimality analysis algorithms over DDs. 

\subsection{Modeling}\label{sec:exact_modeling}

One of the most common purposes of DDs is to model complex combinatorial structures. A DD can represent any combinatorial problem by enumerating solutions as paths. This characteristic is particularly appealing for problems that consider constraints that are usually hard to represent with standard methodologies, e.g., non-linear inequalities. 

A simple procedure to represent a combinatorial problem as a DD is to enumerate all the solutions in a tree and then merge nodes with equivalent paths to the terminal. This procedure is a na\"ive approach that is impractical in many applications due to the exponential growth of the solution set with respect to the number of variables. Thus, most works create DDs using algorithms that avoid explicitly enumerating all the solutions. 

This section discusses general procedures to encode combinatorial problems into  DDs, in particular the recursive formulation strategy presented in Section \ref{sec:pre_recursive}.  We then review modeling techniques that integrate DDs into IP and CP methodologies. Lastly, we present recent extensions of DDs to model problems with continuous variables.

\subsubsection{Recursive Models.}

\cite{hooker2013decision} studied the relationship between DDs and dynamic programming, showing that  DDs can be seen as a compact representation of the state-transition graph (see Section \ref{sec:pre_recursive}). Thus, we can create a DD by building the state-transition graph and merging nodes representing equivalent partial solutions (e.g., the set of paths from the merged nodes to $\terminalnode$ are identical). However, this algorithm can be computationally intractable depending on the size of the state-transition graph.

\cite{bergman2016discrete} revisited this idea and presented a general procedure to create DDs based on top-down algorithms to build relaxed DDs (see Section \ref{sec:approx_construction} for further details). The authors presented recursive models for several classic combinatorial problems (e.g., the maximum independent set) and showed that their procedure can efficiently create exact and relaxed DDs based on such recursive models.  \cite{hooker2017job} analyzed the relationship between DDs and recursive models for sequencing problems. The author formalized some of the ideas introduced by \cite{bergman2016discrete} to create valid DD relaxations and presented a general framework to define recursive models that are suited for exact and relaxed DDs. 

Most papers in this survey use a recursive model to create a DD for two reasons. First, several discrete problems have a natural recursive formulation, e.g., the knapsack problem and many sequencing problems. Second, there exists a wide range of algorithms to construct a DD from a general recursive formulation \citep{bergman2016discrete}. Nonetheless, there are specific cases where other mechanisms are more suited to create a DD encoding, e.g., global constraints that represent a list of feasible solutions \citep{cheng2008maintaining}.

\subsubsection{Network Flow Formulation.} 

One of the most appealing characteristics of DDs for the mathematical programming community is their network flow reformulation \citep{behle2007binary}. This formulation can be integrated into any IP model by adding new variables and constraints.  Moreover, the network flow model is a core component for advanced procedures that combine DDs and IP methodologies, e.g., cutting planes, Benders decomposition, and Lagrangian relaxations. 

Given a DD $\dd=(\nodes,\arcs)$, the network flow model $\flow(\dd)$ uses variables $\by \in \R_+^{|\arcs|}$ to represent the flow traversing each DD arc and the original variables $\bx\in \R^n$ to limit the flow in each layer. Equalities \eqref{eq:flow_all_lit} and \eqref{eq:flow_rt_lit} are balance-of-flow constraints over $\dd$. Constraint \eqref{eq:flow_conv_lit} links the arcs of $\dd$ with solutions $\bx$. 
\begin{subequations}
	\begin{align}
		\flow(\dd)= \{  (\bx; \by) &\in \R^n\times \R_+^{|\arcs|} :  \nonumber \\
		& \sum_{a\in \outgoing(u) } y_a - \sum_{a\in \incoming(u)} y_a = 0, 		&& \; \forall u \in \nodes\setminus\{\rootnode, \terminalnode\}, \label{eq:flow_all_lit} \\ 
		& \sum_{a\in \outgoing(\rootnode) } y_a = \sum_{a\in \incoming(\terminalnode)} y_a = 1, \label{eq:flow_rt_lit}\\
		& \sum_{a \in \arcs: \source(a) \in \nodes_i} \arcValue_{a} y_a = x_i,  	&& \; \forall i \in \varIndex \label{eq:flow_conv_lit} \bigg \},
	\end{align}
\end{subequations}
where $\outgoing(u)$ and $\incoming(u)$ are the outgoing and incoming arcs at a node $u$, respectively. Intuitively, the flow over each path $p\in \paths$ can be seen as the weight of its corresponding point $\bx^p$. By restricting the total flow to have value one, the flow variables represent a convex combination of the points in $\sol(\dd)$. In particular, the polytope $\flow(\dd)$ projected over the $\bx$ variables is equivalent to the convex hull of all solutions represented by $\dd$, i.e., $\projx(\flow(\dd)) = \conv(\feasible)$ \citep{behle2007binary,tjandraatmadja2019target}. Thus, the network flow model $\flow(\dd)$ is an ideal linear formulation of the solution set of $\dd$. In particular, $\flow(\dd)$ can be used to create extended linear formulations of complex combinatorial structures as we review in what follows.

There is a strong relationship between DD network flow models and arc flow formulations based on dynamic programming models \citep{martin1987generating,de2021arc}. As previously mentioned, a DD $\dd$ can be seen as a compact representation of a state-space graph of a recursive model. Thus, $\flow(\dd)$ corresponds to the arc flow formulation for this compact state-space graph. These two lines of research mainly differ on how the graphical structures are constructed and manipulated. For example, DDs have generic reduction algorithms to build exact diagrams with the minimum number of nodes and arcs \citep{bryant1992symbolic}, while graphical manipulations for arc flow formulations are mostly problem specific \citep{de2021arc}. Moreover, there are methodological differences when creating relaxations of the original feasibility set instead of an exact representation (see Section \ref{sec:approx_bounds} for a further discussion).

We now review notable applications using this reformulation as a modeling tool. 
Recent works have shown the advantage of using a DD network flow formulation as a modeling mechanism in a wide variety of real-world applications. \cite{cire2019network} tackled a clinical rotation scheduling problem for medical students, creating a DD-based network flow model to represent all feasible schedules coupled with additional constraints to model other problem characteristics. Another notable real-world application is the design of a heat exchange circuit \citep{ploskas2019heat}. The authors created a DD to represent all possible tube configurations and used its network flow formulation in an MILP model.

While these last two works create a single DD to represent the complex combinatorial component of the problems, several works consider multiple DD network flow models together. \cite{bergman2016decomposition} first proposed the idea of decomposing a problem using multiple DDs that represents specific aspects of the problem. Their procedure creates a network flow model  $\flow(\dd)$ for each DD $\dd$ where the $\bx$ variables are common to each model, i.e., \eqref{eq:flow_conv_lit} are linking constraints that synchronize the solutions among all DDs. \cite{lozano2020consistent} studied the complexity of this multiple network flow model and showed that it is NP-hard in the general case. The authors also proposed a cutting-plane algorithm that solves a maximum flow problem over the DDs to derive cuts and solve the problem more efficiently.

Despite its complexity, the idea of using multiple DDs has been particularly successful in representing non-linear functions as network flow models. \cite{bergman2018nonlinear} employed this procedure to represent non-linear objective functions that admit a recursive formulation. Their technique assumes that the objective function corresponds to the sum of non-linear functions and considers one DD for each function. \cite{bergman2021quadratic} presented a similar decomposition for quadratically constrained problems. Their procedure decomposes the matrix of a quadratic constraint as the sum of multiple smaller matrices, where a DD encodes the solution set induced by each sub-matrix. The authors then used a network flow formulation with linking constraints to linearize the quadratic constraint. Lastly, \cite{nadarajah2020network} created an approximate linear program in the context of DP by employing multiple DD-based network flow models, showing that stronger linear reformulations can be obtained if the (merged) states of the DDs are taken into acount when generating the model.

Recent works also employ network flow models based on multiple DDs to solve challenging applications. \cite{serra2019last} considered a problem that assigns train trips and commuter vehicles to an uncertain set of passengers to minimize the number of commuter trips and total traveling time. Their approach considers a DD for each possible destination and scenario to model the set of passengers associated with a commuter vehicle. \cite{bergman2019binary} addressed the bin packing problem with minimum color fragmentation by building one DD for each bin and creating an IP formulation that links the solutions represented in each DD network flow model. Since all bins are identical, the authors also proposed an IP formulation that uses a single DD network flow model with some modifications to consider multiple bins \citep{mehrani2021models}.
\cite{hosseininasab2019exact} used a similar strategy to tackle the multiple sequence alignment problem with a DD flow model to represent all pairwise sequence alignments and linking constraints to synchronize the DD solutions. The problem is then solved using a Logic-Based Benders Decomposition (LBBD), where the master problem corresponds to the DD flow models with linking constraints and the sub-problems enforce additional constraints over the chosen alignments. Lastly, \cite{lozano2020decision} employed  multiple DDs for the paired job scheduling problem, where each DD represents the specific constraints for a job. The authors proposed a lifted reformulation based on the network flow model where the flow variables are projected out. The main advantage of their lifted reformulation is that it can be constructed without explicitly constructing the DDs.

While the aforementioned works employ the network flow model of \cite{behle2007binary}, other authors have presented variants for DDs that encode stochastic constraints.  \cite{latour2017combining,latour2019stochastic} represented probabilistic constraints with DDs where the parameters follow a probability distribution that depends on the decision variables. The authors used the DD to model the probability of each constraint by encoding the decision variables and the stochastic parameters within the DD. The DD is reformulated into a quadratic constraint model that is linearized and introduced into a MILP formulation of the problem.

\cite{haus2017scenario} also considered a variant of the network flow model for a class of two-stage stochastic programs. Their problem considers endogenous uncertainty, i.e., the first stage decision influence the stochastic process. Their proposed procedure aggregates scenarios with equal cost using multiple DDs and relates these scenarios with the first stage decisions. Similar to \cite{latour2017combining}, each DD computes the probability of achieving a cost value using a MILP reformulation. 

We note that most of the works described here create a DD network flow model to represent a subset of constraints that are generally hard to encode with linear inequalities. For example, several works model the sequencing aspect of the problem  using a DD \citep{cire2019network,ploskas2019heat,serra2019last,hosseininasab2019exact} since other alternatives would involve big-M constraints that have a loose linear relaxation. Alternatively, some authors employ DDs as linearization tools for non-linear expressions \citep{bergman2018nonlinear,bergman2021quadratic} and probabilistic structures \citep{latour2017combining,latour2019stochastic,haus2017scenario}. Thus, a DD network flow formulation is most beneficial when standard procedures lead to poor relaxations and the DD encoding is small.

\subsubsection{Global Constraints.}

In contrast to IP technologies, CP solvers represent combinatorial problems using global constraints, i.e., general-form constraints that encode sub-structures of the problem \citep{rossi2006handbook}. A global constraint has an inference procedure that retrieves information about feasible variable assignments and a propagation mechanism to update the set of feasible solutions inside the constraint. From this perspective, a DD that represents a set of solutions is a global constraint where the inference procedure checks the arc values to determine the current domain of the variables and the propagation procedure removes arcs to update the feasibility set. 

\cite{andersen2007constraint} introduced a general framework for exact and relaxed DDs in CP solvers. Their DD implementation encodes the complete solution set (i.e., the domain store) and propagates the branching decisions. Since the DD could grow exponentially, the authors approximate the solution space using a relaxed DD. Alternatively, we can represent sub-structures of the problem using exact DDs.

One of the main advantages of the DD representation of a global constraint is that it achieves generalized arc consistency (GAC). We say that a variable $x$ is GAC for a constraint $\mathcal{C}$ if for every value of its domain there exists a feasible solution with respect to $\mathcal{C}$. Since the DD represents all feasible solutions of a global constraint $\mathcal{C}$, we can check in polynomial time if a variable assignment is part of a feasible solution or not \citep{andersen2007constraint}. Another important characteristic is that CP solvers construct a DD only once and can efficiently update the graph to eliminate infeasible assignments according to other constraints or branching decisions. We note that these two properties do not hold for relaxed DDs and, thus, it is necessary to build sophisticated procedures to efficiently integrate relaxed DD with CP technologies when exact representations are too large. 

Several researchers have represented different global constraints with DDs and even pre-compile sub-problem structures into DDs to enhance propagation \citep{de2019compiling}. This line of research has inspired new DD variants that are suitable for CP technologies, such as non-deterministic DDs \citep{amilhastre2014compiling}. In the following, we review several works that  represent global constraints using exact DDs. Section \ref{sec:approx_propagation} discusses encoding global constraints with relaxed DDs. 

One of the main applications of DDs in the CP community is to encode global constraints that explicitly enumerate the set of solutions. \cite{cheng2008maintaining} first proposed to model \texttt{Table} constraints (i.e., a list of feasible solutions) using an exact DD. The authors presented an algorithm to convert a \texttt{Table} constraint into a DD by first representing the set of solutions in a tree and then merging identical sub-trees to obtain a reduced DD. Their DD construction procedure was later improved by \cite{cheng2010mdd} and generalized to consider negative \texttt{Table} constraints (i.e., a list of infeasible solutions).

These preliminary works became the stepping stone for DD-based \texttt{Table} constraints and similar structures. \cite{perez2014improving} presented a new algorithm for DD inference over \texttt{Table} constraints that shows superior run-time performance compared to existing methodologies. The same authors also introduced novel algorithms to construct DDs from specialized \texttt{Table} constraints \citep{perez2015relations}. Moreover, the authors presented improved procedures to manipulate DDs, e.g., algorithms to reduce a DD or to add/remove solutions \citep{perez2015efficient,perez2016constructions}, and introduced a parallelization strategy for these algorithms \citep{perez2018parallel}. Lastly, \cite{perez2017soft} studied DD-based propagation for linear cost functions. Their approach improves the cost-based propagator of a DD (i.e., propagation of a linear cost function) and introduces a DD propagator for soft constraints (i.e., constraints that allow infeasible solutions with a cost penalty).

\cite{verhaeghe2018compact} also studied how to efficiently encode \texttt{Table} constraints into DDs. Their work proposes a related data structure called semi-DD (or sDD) where the middle layer is non-deterministic. Intuitively, an sDD is a DD obtained by connecting the leaf nodes of two trees representing partial solutions for half of the variables. The main advantage of this structure is that the maximum number of nodes in each layer can be exponentially smaller than a standard DD. The authors introduced several mechanisms to construct and manipulate an sDD, including a reduction procedure and an algorithm to remove solutions. In a related work, \cite{verhaeghe2019extending} proposed a new DD variant to model \texttt{Smart-Table} constraints. Their basic smart DDs (or bs-DDs) can represent unary constraints (e.g., $x_1\geq 1$) over arcs to compactly encode sets of feasible values. Therefore, bs-DDs avoid having multiple arcs with the same source and target node. \cite{vion2018mdd} also proposed an alternative DD structure for \texttt{Table} constraints. Their work presents a procedure to transform any MDD to special type of BDD where propagation and filtering procedure can be performed more efficiently. In particular, their work adapts filtering procedures proposed for MDDs and other tree structures and shows better results in terms of theoretical complexity bounds and empirical solution times.

Researchers have also employed exact and relaxed DDs to build novel global constraints. \cite{roy2016enforcing} introduced a new global constraint based on DDs that models binary relations over sequences of temporal events. Their empirical results show that their DD-based global constraint is superior to a classical scheduling representation of the same temporal relations. We refer the reader to Section \ref{sec:approx_propagation} for further examples of relaxed DD encodings of global constraints.  

\subsubsection{Continuous Variables.}

So far we have discussed how to model combinatorial problem using DDs. While DDs can represent any bounded set with integer points, until recently there was no DD representation of problems with continuous variables.

\cite{davarnia2021strong} first considered modeling  a mixed-integer set $\feasible\in \R^n$ with  a DD to build outer-approximations.  This work shows that it is possible to construct a DD $\dd$ such that $\conv(\feasible)\subseteq\conv(\sol(\dd))$ where $\dd$ has a finite number of arcs. The main result relies on the fact that it is sufficient to have at most two arcs between consecutive nodes $u$ and $v$ (i.e., $u \in \nodes_i$ and $v\in \nodes_{i+1}$ for some $i \in \varIndex$), where the arc labels take the minimum and maximum feasible value of their associated variable $x$. This type of DD is called an \textit{arc-reduced DD} since it omits arcs with labels other than the minimum and maximum domain values. Therefore, any variable $x\in [a,b]$ can be represented in an arc-reduced DD with at most two arcs emanating from nodes in the layer associated with $x$.

\cite{davarnia2021strong} shows how to build arc-reduced DDs that over-approximate $\conv(\feasible)$ for any $\feasible\in \R^n$, but there are no guarantees on the tightness of such approximations. \cite{salemi2021structure} studied this problem and identified necessary and sufficient conditions to exactly represent $\conv(\feasible)$ with an exact DD. In particular, their paper states that a set  $\feasible\in \R^n$ is \textit{DD-representable} (i.e., there exists a DD $\dd$ such that $ \conv(\feasible)=\conv(\sol(\dd))$) if and only if $\feasible$ admits a  rectangular decomposition (we refer to the original paper for a formal definition). While determining if  $\feasible$ admits a  rectangular decomposition can be challenging, the paper points out that any compact integer  set $\feasible \in \Z^n$ fulfills this property. Moreover, the authors showed that any compact mixed-integer set $\feasible \in \R^n$ with a single continuous variable is DD-representable. This last result is particularly important for decomposition algorithms, as we discuss in Section  \ref{sec:exact_feasible}.

\subsection{Feasibility Checking}\label{sec:exact_feasible}

Another important use of exact DDs is checking the feasibility of a solution. This general idea attracted the attention of researchers in discrete optimization since there are several feasibility checking procedures for IP and CP technologies. The main advantage of using exact DDs for feasibility checking is that we construct the DD only once and repeatedly call it for each new candidate solution.

While feasibility checking is most commonly used in conjunction with IP and CP solvers, two recent works employ this idea without these technologies. In both applications, the solution generation procedure first solves a relatively easier problem to generate solutions fast, and the DD encodes the set of constraints that are hard to represent. \cite{nishino2015bdd} used a DD for feasibility checking to solve the constrained shortest path problem. They employed Dijkstra's algorithm \citep{dijkstra1959note} to solve a shortest path problem but in each step a DD checks the feasibility with respect to the additional constraints. 
The second work generates solutions for a traveling salesman problem (TSP) with preferences using a generative adversarial neural network (GAN) \citep{xue2019embedding}. Since the GAN solutions might violate some of the TSP constraints, the authors used a DD to identify infeasible candidate solutions and guarantee feasibility. 

The following subsections review feasibility checking algorithms in the IP and CP community where DDs are suitable alternatives to other methodologies. As with the procedures described above, the DD can compactly encode the set of hard constraints for the problem and check feasibility in polynomial time with respect to its size. Thus, these techniques are suited for combinatorial structures that lead to small DDs or where alternative procedures (e.g., an IP model) take exponential time to check feasibility.  

\subsubsection{Cutting Planes.}

Cut generation is a mathematical programming procedure to separate infeasible points from the solution set by generating inequalities that are valid for the problem but violated by the infeasible points. This separation problem is NP-hard in the general case, so most cutting plane procedures focus on specific problem structures or use a relaxation of the original problem to generate valid cuts \citep{cornuejols2008valid}. In this context, DDs are an attractive method to solve the separation problem since they provide compact representations of several combinatorial problems, and their network flow formulation lead to cut generation linear programs (CGLP).

We now describe a simple procedure to separate points for a discrete set $\feasible \subseteq \Z^n$ using a DD $\dd$ that exactly represents $\feasible$ \citep{becker2005bdds}. Given a point $\bx'\in \R^n$, we can identify if $\bx' \notin \conv(\feasible)$ if there exists a vector $\blambda \in \R^n$ such that $\blambda ^\top \bx' > \lambda^*$, where $\lambda^* = \max\{\blambda ^\top \bx : \;\; \bx \in \conv(\sol(\dd)) \}$. The main advantage of this separation problem is that computing $\lambda^*$ can be done in polynomial time (w.r.t. to DD size) using a longest-path procedure over the DD with $\blambda$ coefficients. Moreover, the resulting cut $\blambda ^\top \bx \leq \lambda^*$ is valid for every $\bx\in \conv(\feasible)=\conv(\sol(\dd))$.

\cite{becker2005bdds} first studied this DD-based separation problem and proposed a procedure to find vector $\blambda$ that maximizes distance $|\blambda^\top \bx' - \lambda^* |$ for any $\bx'\notin \conv(\feasible)$. The authors used an iterative sub-gradient procedure that computes $\lambda^*$ using a longest-path over $\dd$ and updates $\blambda$ with the longest-path solution. \cite{davarnia2020outer} studied this procedure to create outer-approximations for non-linear problems and enhanced the sub-gradient routine by adding normalization constraints to improve convergence. 

This simple DD-based separation algorithm can also be reformulated as a CGLP using the dual of the network flow model $\flow(\dd)$  \citep{behle2007binary}. While this CGLP might be computationally slower than a sub-gradient routine \citep{davarnia2020outer}, it is easier to implement and can be modified to obtained structural properties for the resulting cuts. Specifically, two recent works develop novel DD-based CGLPs based on \cite{behle2007binary} CGLP. \cite{tjandraatmadja2019target} proposed a DD-based CGLP to generate target cuts based on the polar set of $\flow(\dd)$. Their CGLP can compute facet-defining inequalities by just considering a few additional variables and constraints. In contrast, \cite{castro2021combinatorial} introduced a reformulation of $\flow(\dd)$ for binary sets and created a CGLP that resembles a maximum flow problem. The authors also proposed a combinatorial approach to generate cuts that solves a max-flow problem over the DD. Their procedure is orders of magnitude faster than DD-based CGLPs and the sub-gradient routines but might not cut all infeasible fractional points $\bx' \notin \conv(\sol(\dd))$.

In addition to the previously mentioned approaches, \cite{behle2007binary} proposed two types of DD-based logical cuts for branch-and-cut routines. The first one returns an exclusion cut when the DD is infeasible given the current branching decisions (i.e., the DD has no \rt\ paths/feasible solutions due to branching decisions). The second one is an implication cut that is returned when the values of some variables are fixed after updating the DD with the branching decisions.

Cutting-plane procedures are usually coupled with constraint enhancement routines so that the resulting cut removes a larger portion of the fractional space. However, there is currently a limited number of DD-based algorithms to enhance valid inequalities. \cite{becker2005bdds} introduced a simple procedure to modify the coefficients of a valid inequality to increase its face dimension. However, the resulting inequality might not separate the set of points that the original constraint does. \cite{behle2007binary} proposed a DD-based lifting procedure for cover cut inequalities based on classic techniques \citep{wolsey1999integer} in which the coefficient sub-problem is solved using a DD instead of an integer or linear program. 
Recently,  \cite{castro2021combinatorial} introduced a general lifting procedure for combinatorial problems where a DD iteratively tilts valid inequalities in order to increase their dimension. Their procedure resembles the simple approach by \cite{becker2005bdds} but the lifted inequality is guaranteed to dominate the starting inequality. 

As with other DD procedures reviewed so far, these cutting plane techniques leverage the fact that the DD is built only once and can be used to find optimal solutions with different objective functions. For example, the sub-gradient routines use this property to create new sub-gradients based on the DD optimal solutions. The same is true for the CGLP models and the inequality enhancement procedures. We note that these cutting plane techniques are also valid if we replace the exact DD with a relaxed DD \citep{tjandraatmadja2019target,castro2021combinatorial,davarnia2020outer}. In this case, the generated inequalities are valid for the original problem but the cutting plane procedure will not separate all infeasible points.

\subsubsection{Benders Decomposition.}

Benders decomposition \citep{benders1962partitioning} is another IP methodology  that has benefited from DDs. The overall idea is to decompose the problem into a master problem and one or more sub-problems, where each sub-problem is represented by a DD. Each DD $\dd$ can then be used to generate logic-based cuts or traditional Benders cuts with the help of the network flow model $\flow(\dd)$.
\cite{bergman2021quadratic} first studied this approach for quadratically constrained integer problems. Their procedure decomposes the quadratic matrix into multiple smaller matrices where a DD can easily represent the solution set of each component matrix. Their work proposes a Benders decomposition scheme where each sub-problem solves the shortest path problem over its corresponding DD to generate a Benders cut.

Similar ideas have also been used to solve stochastic programming problems.  \cite{lozano2018binary} introduced a DD-based Benders decomposition for a family of two-stage binary stochastic programming problems and applied them to a stochastic TSP variant. Their approach uses a DD to represent the sub-problem (i.e., one for each scenario) and solves a max-flow problem over the DD to create Benders cuts. \cite{guo2021logic} applied this same procedure to tackle a stochastic distributed operation room scheduling problem. 


In these applications the DD represents a combinatorial problem that might take a considerable time for an IP model to solve. Since the sub-problems of these Benders decompositions are relatively small (e.g., 20 to 40 cities for the traveling salesman problem), the resulting DD is small enough to store in memory and can solve the sub-problems in fractions of a second. Since the set of solutions in the DD remains the same, the overhead of constructing the DD is negligible if we consider that we need to solve the sub-problem thousands of times. Moreover, the DD network flow model is a convex reformulation of the original discrete non-convex sub-problem. Thus, this model can be used to generate LP-based Benders cuts, which is impossible for the original discrete sub-problem.

Alternatively, \cite{salemi2021structure} proposed a DD-based Benders decomposition where the master problem is given by a DD and the sub-problems are, for example, linear programs. Their procedure relies on the fact that mixed-integer programming models with a single continuous variable are DD-representable. Thus, we can create a DD for the  master problem with the integer variables of the original problem and a single continuous variable to approximate the cost of the sub-problems. Their decomposition, however, needs to update the master problem DD for each sub-problem cut, which can be computationally expensive. To mitigate the potential exponential growth of the exact DD, the authors proposed an iterative DD-based Benders approach that uses relaxed and restricted DDs. 

\subsubsection{Inference.}

The CP community has also considered DDs for feasibility checking to infer no-good assignments, i.e., an infeasible set of variable assignments \citep{schiex1994nogood}. CP solvers employ no-goods to avoid exploring portions of the solution space that it has already shown to be infeasible.  

\cite{subbarayan2008efficient} first proposed the idea of using DDs to infer no-good assignments in CP solvers. The author showed that the problem of finding a minimum no-good over a DD is NP-hard and proposed a heuristic procedure that traverses the DD to find small sets of no-good assignments. \cite{gange2011mdd} revisited this idea and presented an incremental algorithm to find no-goods that searches just a portion of the DD in the vicinity of the last set of arcs removed from a DD. \cite{gange2013explaining} extended this no-good algorithm for cost-based reasoning, i.e., to identify assignments that lead to sub-optimal solutions. A similar idea was implemented for SAT solvers \citep{kell2015clause} where a  DD represents a subset of the constraints (i.e., clauses) and identifies no-goods for clause generation. 

We note that no-good inference is related to the cutting plane procedures in IP technologies. A no-good set can be encoded as a linear constraint and added into an IP solver in a branch-and-cut procedure. In fact, the DD-based logic constraints introduce by \cite{behle2007binary} are a particular type of no-good that a DD infers during the search. Thus, these advanced technologies for no-good inference in the CP literature can also benefit IP solvers. 

Alternatively, DDs can be used to efficiently check if one or more solutions satisfy a global constraint. \cite{jung2021checking} proposed an inference procedure for this task that relies on a novel DD inclusion operator. Their approach creates two DDs, one for the solution set and one for the constraint, and then employs their inclusion operator to check if the set of solutions satisfies the constraint. The authors showed that their inclusion operator is  more efficient than a DD intersection operator when applied to some well-know global constraints, e.g., the \texttt{Sequence} constraint and \texttt{GCC} (global cardinality constraint). 

\subsection{Solution Extraction} \label{sec:exact_extract}

In contrast to feasibility checking, we can use a DD to generate candidate solutions for a given problem. Since DDs represent all the solutions of a combinatorial problem, we can easily extract solutions by choosing any path in the DD. In particular, we can use this property to extract solutions with certain structure or solutions that are optimal for a specific linear objective function. 

\cite{hadzic2004fast} first explored this idea for a manufacturing problem. The authors created a DD to represent all product configurations and proposed a polynomial-time algorithm to extract solutions that satisfy a set of configuration requirements specified by a user. Their algorithm ignores all arcs that represent infeasible configurations and returns a sub-diagram without the ignored arcs. While the authors tested their procedure in a manufacturing problem, the same idea can be applied for any other combinatorial problem to extract solutions with user-specified variable assignments. 

Solution extraction has the advantage that we can decompose the problem into two parts. The first part, encoded with a DD, represents a relaxation of the original problem to create candidate solutions. The second part of the problem checks the feasibility of the extracted solution and gives feedback to the DD to extract new solutions. This mechanism is a generalization of the column generation procedure and, thus, it can be applied with other technologies. We also point out that this decomposition scheme is suited for problems where the solution extraction and feasibility checking are relatively easy to solve and a DD compactly represents the set of solutions. 

\subsubsection{Column Generation.}

Solution extraction has been mostly used as a sub-routine of the column generation procedures in IP/LP technologies. \cite{morrison2016solving} proposed a general scheme that applies DDs in a branch-and-price algorithm \citep{barnhart1998branch}. The DD solves the pricing problem of the column generation sub-routine, i.e., the DD returns a promising solution and adds a new column to the master problem.  The algorithm then separates the last queried solution from the DD to avoid duplicated solutions. Their implementation also considers a branching rule over the original variables, which can be easily integrated into the DD by ignoring arcs with specific values. 

There are three advantages of the procedure by \cite{morrison2016solving}. First, the DD can solve the pricing problem in polynomial time (with respect to its size) for any linear objective function. Thus, their implementation constructs the DD only once and can solve the pricing problem multiple times. Second, the DD can return all optimal solutions for the pricing problem and, thus, add multiple columns in each iteration.  Third, branching decisions can be easily included to the DD without changing the complexity of the pricing problem. Therefore, this procedure can be a better alternative to other methodologies when there is a compact DD representation of the pricing problem. 

\cite{morrison2016solving} tested their procedure on the graph coloring problem. \cite{kowalczyk2018branch} implemented the DD-based branch-and-price procedure for a parallel machine scheduling problems. \cite{raghunathan2018seamless} applied a similar branch-and-price approach to solve a transportation scheduling problem. Their application has multiple pricing problems, each one modeling the tours of an origin-destination pair using a DD. Lastly, \cite{riascos2020branch} employed the DD-based branch-and-price technique to solve a kidney exchange problem where multiple DDs represent the set of possible chain donations of different sizes.

\subsection{Solution-Space Analysis}

DDs provide a compact representation of the solution space, which is also useful if we want to enumerate and analyze the set of solutions. For example, we can analyze the set of (near) optimal solutions or study the convex hull of the solutions of a DD. We now review procedures to analyze the solutions encoded in a DD and to enumerate solutions with specific characteristics.

\subsubsection{Post Optimality Analysis.}

\cite{hadzic2006postoptimality} studied this topic for discrete optimization problems and  introduced three procedures over an exact DD. The first one is a cost-based domain analysis that identifies the set of variable assignments that are present in at least one near-optimal solution.  The second is a conditional cost-based domain analysis, which restricts a subset of variables to take specific values and then performs cost-based domain analysis over the subset of solutions. Lastly, their frequency analysis computes the percentage of solutions with a particular variable assignment. 

Since representing the set of solutions using a DD is intractable for most combinatorial problems, \cite{hadvzic2007cost} introduced the concept of sound DDs for post-optimality analysis. A DD is sound for a specific post-optimality analysis if it yields the same results that an exact DD would. Thus, a sound DD might represent a larger solution set than an exact DD but preserves certain properties to correctly perform the analysis. The authors presented a pruning and contraction procedure to create sound DDs for cost-based domain analysis. While the resulting sound DD is significantly smaller than an exact DD, their procedure requires a starting DD that is either exact or represents all feasible solutions within a cost range. Thus, creating sound DDs with this procedure is intractable for large problems.

Recently, \cite{serra2019compact} revisited the idea of sound DDs for post-optimality analysis and provided new insights. The authors presented several theoretical results related to sound DDs, including a sound-reduction algorithm that constructs the smallest sound DD. They also introduced a general construction procedure for ILP models that creates a sound DD from a branching tree. Their experiments over a wide range of ILP benchmarks show that the resulting sound DD is a more suitable alternative to post-optimality analysis than, for example, branch-and-bound enumeration. 

\subsubsection{Solution Enumeration.}

Recent works have also employed DDs to represent the set of non-dominated solutions for a  multi-objective discrete optimization problem, i.e., the Pareto frontier. \cite{bergman2016multiobjective} first tackled the problem by representing the set of feasible solutions with an exact DD and then enumerating the non-dominated solutions using a multi-criteria shortest path algorithm over the DD. This work was extended by \cite{bergman2021network} where the authors presented three procedures that modify the DD while preserving the set of non-dominated solutions. The authors also proposed a bidirectional multi-criteria shortest path procedure to enumerate the non-dominated solutions, which outperforms the unidirectional approach \citep{bergman2016multiobjective}. \cite{suzuki2018fast} presented a similar procedure to enumerate non-dominated solutions for the 0-1 multi-objective knapsack problem. Their algorithm encodes all feasible solutions using a ZDD and employs specialized filtering procedures based on the knapsack structure to prune dominated solutions.

All the papers discussed so far enumerate solutions to analyze optimal or near-optimal solutions. Conversely, we can employ the DD structure to obtain insightful information on the combinatorial set and, thus, ignore the objective function.  \cite{suzuki2018exact} built a DD to compute the strongly connected reliability of a network, i.e., the probability that the network will remain strongly connected when removing one or more edges. The authors showed how to create a DD that enumerates all possible strongly connected sub-networks and how to compute the desired probability by traversing the DD once. Thus, this work employed DDs to get structural properties of the connectivity problem over a network. \cite{haus2017compact} constructed a DD to encode all the members of an independent set to analyze the size of the solution set. The authors showed that the DD representation has polynomial size in the number of variables for packing and set covering problems with specific characteristics. Thus, this work opens the possibility of solution set analysis for problems with tractable DD size.

\subsubsection{Polyhedral Analysis.}

We now review two works that employ DDs for polyhedral analysis, an important area of research in mathematical programming due to its relevance for solving IP problems. \cite{behle2007facet} introduced a procedure to enumerate vertices and facets of the convex hull of a DD solution set, i.e., $\conv(\sol(\dd))$. Their technique considers a binary solution set in an exact DD, where each path in the DD corresponds to a 0/1 vertex of the convex hull polytope. The enumeration procedure starts with an initial facet that is rotated over the DD to obtain a new facet.

\cite{tjandraatmadja2019target} also presented a polyhedral analysis procedure based on DDs to certify the dimension of any inequality that is a face of $\conv(\sol(\dd))$. Their procedure finds a set of affine independent points in the DD that are tight for the face by solving a flow problem over the DD. Their procedure then uses a heuristic approach to generate affine independent points based on the flow values of each arc in the DD. The number of affine independent points gives a lower bound on the dimension of the face. 

We note that polyhedron analysis is intractable for general IP models since constructing the convex hull of all the solutions is NP-hard \citep{wolsey1999integer}. DDs provide a more manageable procedure to enumerate all the solutions and, thus, give valuable insights to the polyhedral structure of problems that have a compact DD encoding.

\section{Approximate DDs} \label{sec:approxDD}

Relaxed and restricted DDs are limited size DDs that over- and under-approximate the set of feasible solutions, respectively (see Section \ref{sec:pre_approx} for a formal definition). These graphical structures are desirable when tackling large-scale problems since the size of an exact DD can grow exponentially with the number of variables. Specifically, approximate DDs provide primal and dual bounds for a problem and, thus, can be embedded in a search procedure to find optimal solutions for large problems. 

This section presents an overview of existing procedures to construct and employ DD approximations. We distinguish three main research areas related to approximate DDs (see Table \ref{tab:lit_approx}). The first research area focuses on DD construction procedures to create tight approximations.  The second area studies bound computation based on DD approximations and specialized search procedures that employ these bounds. The last research area focuses on propagation procedures for CP solvers based on relaxed DDs.

\subsection{Approximate DD Compilation}\label{sec:approx_construction}

While several works have shown the advantage of using approximate DDs \citep{Bergman2016book}, the question of how to construct one that provides tight bounds remains open. The size of the DD plays an important role in the quality of the approximation, i.e., bigger diagrams are expected to produce tighter bounds. Nonetheless, two DDs with equal size can provide significantly different bounds. 

In the following, we introduce the most commonly used procedures to create approximate DDs and discuss recent DD construction algorithms. These procedures can be applied to a large range of problems, in particular to problems that can be reformulated recursively. Also, note that there is no clear superiority between these construction algorithms but their advantages and disadvantages can be used as a guideline to choose the most appropriate procedure for a specific application. 

\subsubsection{Top-down and Iterative Refinement.}

The most common techniques to construct approximate DDs are the top-down and the iterative refinement procedures \citep{Bergman2016book}. The main advantage of both procedures is that they  can be applied to any combinatorial problem that has a recursive formulation \ref{model:gen_rec}. These procedures create an approximate DD $\dd$ of limited size by restricting the maximum number of nodes per layer, i.e., limiting its width $\width(\dd) = \max_{i \in \varIndex}\{|\nodes_{i}|\}$. Thus, we can control the size of the DD approximation by changing the maximum width limit $\maxwidth \geq 0$.

\begin{algorithm}[htb]
	\caption{DD Top-Down Construction} \label{alg:dd_topdown}
	\begin{algorithmic}[1]
		\Procedure{\algTopDownDD}{\ref{model:gen_rec}, $\maxwidth$}
		\State Create DD $\dd=(\nodes, \arcs)$ with $n+1$ empty layers
		\State Create the root and terminal node, i.e.,  $\rootnode\in \nodes_{1}$ and $\terminalnode \in \nodes_{n+1}$ 
		\State Assign the initial state to the root node, i.e., $\bS(\rootnode) = \bS_1$
		\For{ $i \in \{1,\dots, n-1\}$}
		\For{ $u \in \nodes_i$}
		\State Create an arc $a$ emanating from $u$ for each possible value in $\feasibleDP_i(\bS(u))$
		\For{ $a\in \outgoing(u)$ }
		\If{there exists node $u'\in \nodes_{i+1}$ with $\bS(u') = \transitionF_i(\bS(u), \arcValue_a)$}
		\State Direct arc $a$ to node $u'$, i.e., $\target(a) = u'$
		\Else
		\State Create  $u'$ in $\nodes_{i+1}$ with $\bS(u') = \transitionF_i(\bS(u), \arcValue_a)$ and point arc $a$ to $u'$ 
		\EndIf
		\EndFor
		\EndFor
		\If {$|\nodes_{i}| > \maxwidth$} 
		\State \algDDMerge($\nodes_{i}$) (relaxed DD) or \algDDDiscard($\nodes_{i}$) (restricted DD)
		\EndIf
		\EndFor
		\For{ $u \in \nodes_{n}$}
		\State Create an arc $a$ emanating from $u$ for each value in $\feasibleDP_n(\bS(u))$ with $\target(a) =\terminalnode$
		\EndFor
		\State \algReturn\ $\dd$
		\EndProcedure
	\end{algorithmic} 
\end{algorithm}

We first describe the top-down construction procedure \citep{bergman2011manipulating}, which can create exact, relaxed, and restricted DDs by making small changes in the implementation (see Algorithm \ref{alg:dd_topdown}). The algorithms starts with a DD structure with two nodes, the root node and terminal node (lines 2-4). The procedure traverses the layers creating the emanating arcs of each node and the nodes in the following layer (lines 5-12). Note that the procedure creates one node for each reachable state given by \ref{model:gen_rec}, thus, it will return an exact DD if the maximum width $\maxwidth$ is large enough. If the number of nodes in a layer is larger than the maximum width, the procedure will reduce the number of nodes to create either a relaxed or restricted DD (line 13-14). 

In the case of a relaxed DD, the top-down procedure will merge the set of nodes into at most $\maxwidth$ nodes to enforce the DD width limit. The \algDDMerge($\nodes_{i}$) procedure utilizes an appropriate merging operator $\merge$ to guarantee that no feasible solutions are lost and, thus, to construct a valid relaxation for the problem (see Section \ref{sec:pre_approx} for further details on merging operators). In contrast, the \algDDDiscard($\nodes_{i}$) procedure selects a subset of nodes to eliminate from the diagram. By doing so, the resulting restricted DD represents a sub-set of the feasible solutions of the original problem.

The decisions on which nodes to merge/discard are usually done heuristically based on the state information of the nodes. Since these decisions are crucial to create stronger approximations, several researchers have developed general heuristics that lead to better approximations than a random node selection. We note that there are several papers on merging procedures for relaxed DDs but that is not the case for restricted DDs and discarding heuristics \citep{bergman2014heuristic}. This discrepancy is due to the wider exploration of relaxed DDs in the literature compared to restricted DDs.

For the relaxed DD case, \cite{bergman2011manipulating} proposed a simple strategy that merges nodes with respect to their longest/shortest path from the root. Thus, the heuristic avoids merging nodes that can potentially impact the DD bound. A similar idea was then proposed for discarding nodes to avoid removing optimal solutions \citep{bergman2014heuristic}. \cite{frohner2019towards} further studied the \cite{bergman2011manipulating} merging heuristic and proposed tie-breaking rules to merge nodes that have similar relaxed states. \cite{frohner2019merging} also presented a binary classifier procedure that chooses a merging heuristic in each DD layer. Their heuristic outperforms the simple heuristics in terms of dual bounds but it can be harder to implement due to the time needed to train the classifier. Recently, \cite{de2020single} proposed a novel merging procedure for a single-machine scheduling problem that guarantees building an $\epsilon$-approximation of the problem. Their procedure groups states with values in a certain range that depends on $\epsilon$. While the \cite{de2020single} merging procedure was designed for a specific scheduling problem, similar ideas could be applied to other problems.

The main advantage of the top-down algorithm is its flexibility and simplicity. The procedure can be used to construct any type of DD (exact, relaxed, or restricted) given a recursive formulation of the problem. Moreover, the algorithm has shown to be computationally efficient in practice, allowing researchers to construct DDs with large widths in fractions of a second \citep{Bergman2016book}. However, one of it main disadvantages its the looseness of the approximation that it produces. The procedure merges/discards nodes considering only the state information of the current nodes, which can lead to poor approximations and, thus, weak bounds. While there exist alternatives that include a look-ahead step, these procedures can be computationally expensive \citep{horn2021based}.

An alternative is to construct a relaxed DD using the iterative refinement procedure \citep{hadzic2008approximate}. In contrast to the top-down algorithm, this procedure starts with an initial relaxed DD and iteratively increase its width by splitting nodes. The main advantage is that we can use information on the resulting DD to guide the refinement in the next iteration. This procedure can create exact and relaxed DDs but it cannot construct restricted DDs since it only removes infeasible solutions (Algorithm  \ref{alg:dd_iterative}, line 7). 

\begin{algorithm}[htb]
	\caption{DD Iterative Refinement Construction Procedure} \label{alg:dd_iterative}
	\begin{algorithmic}[1]
		\Procedure{\algDD}{\ref{model:gen_rec}, $\maxwidth$, $\merge$}
		\State Create a width-one DD $\dd$
		\While{$\dd$ has been modified}
		\For{$i \in \varIndex$}
		\State Update (relaxed) state information in each node $u \in \nodes_i$ 
		\State Split nodes in $\nodes_i$ and update their relaxed states if needed 
		\State Check all outgoing arcs of layer $\nodes_i$  and eliminate infeasible arcs 
		\EndFor
		\State  Update bottom-up state information in every node $u \in \nodes$  \hfill 
		\EndWhile
		\State \algReturn\ $\dd$ 
		\EndProcedure
	\end{algorithmic} 
\end{algorithm}

Algorithm \ref{alg:dd_iterative} illustrates the iterative refinement procedure given a recursive model \ref{model:gen_rec}, a maximum width $\maxwidth$, and a merging operator $\merge$. The procedure starts by constructing a width-one DD, i.e., a DD $\dd$ where each layer $\nodes_i$, with $i \in \varIndex$, has a single node and emanating arcs for each value in their corresponding domain (line 1). The root node is associated with the initial state (i.e., $\bS(\rootnode) =\bS_1$) and each following node corresponds to a relaxed state computed using the merging operator and its incoming arcs:
\begin{equation} \label{eq:update_top}
	\bS(u) = \bigmerge_{a\in \incoming(u)} \transitionF_{i-1}( \bS(\source(a)), \arcValue_a), \qquad\qquad  \forall u \in \nodes_i,\; i \in \{2,\dots,n+1\}.
\end{equation}

The procedure iteratively refines $\dd$ one layer at a time starting with the first layer until it cannot be updated any further (lines 3-8). It first updates the relaxed states of each node in the current layer using \eqref{eq:update_top}, to guarantee up-to-date state information for each node (line 5). The algorithm uses the relaxed states to split a node $u\in \nodes_i$, i.e., it creates a new node $u'$ with the same outgoing arcs that $u$ has and redirects a portion of the incoming arcs of $u$ to $u'$. The states of the split nodes are then updated and used to identify if any of the outgoing arcs lead to infeasible assignments.

The main advantage of this procedure is that in each iteration we obtain a DD that better approximates the original problem. Moreover, we can use the current DD relaxation to decide how to split nodes in order to remove as many infeasible arcs (i.e., arcs representing infeasible assignments) as possible. We can also create additional relaxed states for a node $u$ using the partial assignments from $u$ to node $\terminalnode$ (i.e., bottom-up relaxed states). These relaxed states are commonly used to find infeasible arcs that are not necessarily identifiable with standard relaxed states \citep{hadvzic2007cost,cire2013multivalued}.

\smallskip
\begin{example}\label{exa:knap_relaxeddd}
	Consider the knapsack problem with feasible set $\feasible=\{ \bx\in \B^4: 7x_1 + 5x_2 + 4x_3 +x_4 \leq 8 \}$ and recursive model \ref{model:knap_rec} from our previous examples. We construct a relaxed DD for $\feasible$ as follows. For each node $u \in \nodes_i$, we consider a relaxed states $\bS(u)=(\Qmin(u), \Qmax(u))$ where $\Qmin(u)$ and $\Qmax(u)$ represent the minimum and maximum load of the knapsack at node $u$ and stage $i \in \{1,\dots,5\}$, respectively. Thus, the merging operator is given by $\merge =(\min, \max)$. We choose this relaxed state representation to identify if a node encodes multiple states and, therefore, if it is a candidate for splitting. Intuitively, a node $u\in \nodes$ represents a single state of \ref{model:knap_rec} if $\Qmin(u) = \Qmax(u)$.
	
	Lastly, we identify if an arc $a\in \arcs$ with source $\source(a)\in \nodes_i$  is infeasible if
	\begin{equation}
		\Qmin(\source(a)) + w_i \arcValue_a >8 \tag{KP-R1}\label{fil:knap}
	\end{equation}
	\noindent for any $i \in \varIndex$. If arc $a$ satisfies condition \ref{fil:knap}, then all paths traversing $a$ represent solutions with a knapsack load above its limit. Thus, we can remove arc $a$ from $\dd$.
	
	\begin{figure}[htb]
		\includegraphics{figs-bdd-vars.pdf}
		\includegraphics{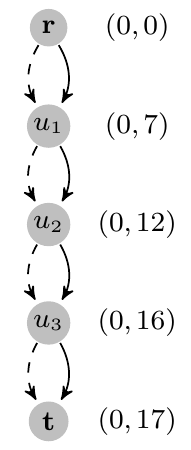}
		\hfill
		\includegraphics{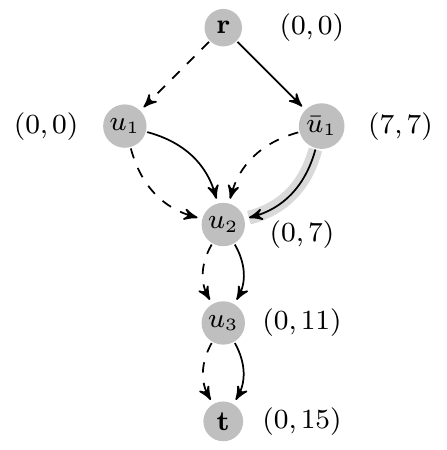}
		\hfill
		\includegraphics{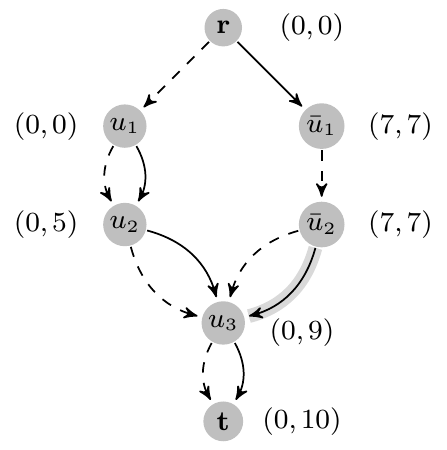}
		\hspace{-2em}
		\includegraphics{figs-bdd-arc-values.pdf}
		\caption{Iterative refinement procedure for $\feasible=\{ \bx\in \B^4: 7x_1 + 5x_2 + 4x_3 +x_4 \leq 8 \}$ and $\maxwidth=2$. The figure shows a width-one DD (left), a DD after splitting layer $\nodes_2$ (middle), and a DD after splitting layer $\nodes_3$ (right). Values next to nodes represent relaxed states. Highlighted arrows correspond to infeasible arcs that can be removed.} 
		\label{fig:knap_dd_const}
	\end{figure}
	
	Figure \ref{fig:knap_dd_const} illustrates the iterative refinement procedure for \ref{model:knap_rec} with maximum width $\maxwidth=2$ and relaxed states as defined above. The left most DD corresponds to a width-one DD, where the values next to each node represent its relaxed state. The middle DD depicts the update, split, and filter sub-routines in layer $\nodes_2$. Node $u_1 \in \nodes_2$ is split into two (i.e., $u_1$ and $\bar{u}_1$), each node with one incoming arc and relaxed states updated accordingly. Notice that arc $a=(\bar{u}_1, u_2)$ with value $\arcValue_a = 1$ is infeasible, since it satisfies \ref{fil:knap}. 
	The right graph shows the resulting DD after refining layer $\nodes_3$. In this DD, condition \ref{fil:knap} identifies infeasible  arc $a=(\bar{u_2},u_3)$ with value $\arcValue_a = 1$.

	While condition \ref{fil:knap} identifies several infeasible arcs during the DD construction procedure, there exist paths in the right most DD of Figure \ref{fig:knap_dd_const} that correspond to infeasible solutions, e.g.,  path $(\rootnode, u_1, u_2, u_3, \terminalnode)$ associated with point $\bx=(0,1,1,1)$ is infeasible. This path cannot be removed from the DD since all its arcs are also associated with feasible solutions, i.e., all the arcs  violate \ref{fil:knap}. We can remove this path from the DD if we further split node $u_2$. \hfill $\blacksquare$
\end{example}

\medskip
Analogous to the top-down procedure, the tightness of a relaxed DD from an iterative algorithm depends on the splitting heuristic. Ideally, we want to split nodes such that the relaxed states are close to exact states. However, an appropriate splitting procedure will greatly depend on the application at hand. For example, splitting nodes to exactly represent specific problem characteristics is a good strategy to identify infeasible assignments and improve the DD relaxation \citep{cire2013multivalued,castro2019mPDTSP,castro2020SolvingDF}. By doing so, the relaxed DD will satisfy a subset of the constraints and, thus, have theoretical guarantees on the quality of the relaxation.

These two construction procedures have key distinctions that make them appealing for different purposes. Top-down procedures are considerably faster to implement and deploy than iterative refinement procedures, so a top-down construction is usually preferable for prototyping new ideas or when the user needs to construct multiple DDs. In contrast, iterative refinement procedures can lead to stronger bounds or exact enforcement of a subset of the constraints, thus, guaranteeing certain properties for the resulting relaxation. However, the iterative procedure is usually computationally slower than the top-down approach due to the multiple refinement iterations, and thus, it is preferable when we need to construct few relaxed DDs.

\subsubsection{Other Construction Algorithms.}

A less commonly used technique to construct relaxed DDs is the separation procedure \citep{cire2014separation}. This algorithm is a variant of the iterative refinement since it starts with a width-one DD and iteratively splits nodes. However, the separation procedure splits nodes systematically so that all the paths in the DD satisfy a constraint or have a longest-path value larger than a certain threshold \citep{bergman2011manipulating}. This technique can generate highly-accurate relaxations if, for instance, we separate the constraint that is currently violated by the longest path in a maximization problem \citep{bergman2016theoretical}. Recent work by \cite{van2020graph,van2021graph} shows that separation procedures are quite effective for the graph coloring problem: the DD relaxation yields bounds competitive with state-of-the-art methodologies.

Besides these construction techniques, two recent works have developed novel procedures to construct relaxed DDs. \cite{romer2018local} proposed a local search framework with operations to merge nodes, split nodes, and  redirect arcs. Their procedure generalizes the top-down and iterative refinement algorithms and consistently yields better bounds than these alternatives. \cite{horn2021based} introduced an \astar-inspired algorithm to construct relaxed DDs. Starting at the root node, the procedure keeps a priority list of nodes to create next and can create and merge nodes in different layers. Recent works showed that the \astar-inspired approach can return tighter bounds than the top-down algorithm for the prize-collecting scheduling problem \citep{horn2021based} and the longest common subsequence problem \citep{horn2021common}.

Lastly, \cite{bergman2017finding} studied the problem of finding the best DD relaxation for a given width limit. The authors presented an IP model that partitions the nodes in each layer to obtain the strongest dual bound. Their experiments show that the resulting DD provides significantly stronger bounds than any other methodology. However, the computational time to solve their IP model can be quite long and, thus, impractical for many applications.

\subsubsection{Variable Ordering.}

One of the main challenges when constructing a DD is the variable ordering, i.e., the assignment of variables to layers. The size of a DD depends on the variable ordering and the optimal ordering can lead to significantly smaller DDs. The problem of finding the optimal variable ordering is NP-hard \citep{bollig1996improving} and has been extensively studied for exact DDs. Several authors propose heuristic variable orderings for different applications, including the knapsack \citep{behle2008threshold} and the maximum independent set problem \citep{bergman2011manipulating}. 

\cite{cappart2019improving} introduced a general ordering procedure based on reinforcement learning (RL). While their technique can be applied to any optimization problem, it is only compatible with the top-down construction procedure: during construction the RL agent chooses the variable that will be assigned to the next layer. A recent work by \cite{parjadis2021improving} showed that the RL-based variable ordering can significantly improve the bound quality and the number of nodes explored during a branch-and-bound search.  

Alternatively, \cite{karahalios2021variable} studied portfolio algorithms to select an appropriate variable ordering for a problem. Their procedure considers several variable ordering heuristics and selects one for a given problem. In contrast to the \cite{cappart2019improving} approach, the portfolio methodology provides a complete variable ordering and can be used with any DD construction algorithm.  

\subsection{DD-based Bounds} \label{sec:approx_bounds}

DD approximations are commonly used to compute bounds for optimization problems. The main advantage of DD bounds is that their tightness can be control by either changing the size of DD or the mechanism to build the diagram (see Section \ref{sec:approx_construction} for more details). Therefore, DD approximations provide a flexibility that is harder to obtained with other procedures (e.g., linear programming or standard primal heuristics).

In the following, we describe the main mechanisms to compute primal and dual bounds with DD approximations and how to enhance these bounds. We then discuss
different ways to embed DD bounds into a branch-and-bound algorithm. 

\subsubsection{Dual Bounds.} \label{sec:approx_dual}

\cite{andersen2007constraint} first introduced the idea of using relaxed DDs to obtain dual bounds for an optimization problem. A relaxed DD over-approximates the set of feasible solutions, so optimizing over the DD provides a valid dual bound for the problem. The authors consider a DD $\dd=(\nodes, \arcs)$  and the optimization problem $\max\{f(\bx): \bx\in \sol(\dd) \}$ where the objective function is separable, i.e., $f(\bx) = \sum_{i \in \varIndex}f_i(x_i)$. The main feature of this problem is that it can be solved using a longest-path algorithm over $\dd$. Intuitively, we assign  each arc $a\in \arcs$ a length given by $\arcLength_a = f_i(\arcValue_a)$ for any arc $a$ emanating from node $u \in \nodes_i$. Then, the longest \rt\ path is the optimal solution over $\sol(\dd)$. 

This simple longest-path procedure provides a valid dual bound for any optimization problem where the feasible region is over-approximated by $\dd$. Note that the longest-path optimization holds for separable objective functions (e.g., a linear objective), and it can be generalized for problems with a recursive formulation \ref{model:gen_rec} \citep{bergman2016discrete,bergman2018nonlinear}. In this case, the length of an arc $a\in \arcs$ corresponds to the immediate reward over $\arcValue_a$ and the relaxed state of its emanating node, i.e., $\arcLength_a = \costF_i(\bS(u), \arcValue_a)$ for any arc $a$ emanating from node $u \in \nodes_i$. The procedure will return a valid dual bound if the node merging operator satisfies the conditions described in Section \ref{sec:pre_approx} for a proper relaxation \citep{hooker2017job}.

Researchers have tested these bounds in a wide variety of combinatorial problems, including set covering  \citep{bergman2011manipulating}, multidimensional bin packing \citep{kell2013mdd}, maximum independent set \citep{bergman2014optimization,bergman2016discrete}, maximum cut, maximum 2-satisfiability  \citep{bergman2016discrete}, and  graph coloring \citep{van2020graph,van2021graph}. These papers create a single relaxed DD to approximate the feasible set and compute bounds using the shortest/longest path procedure over the DD or, alternatively, the DD network flow model. 

One of the most popular applications of relaxed DDs is for sequencing problems. These problems are generally formulated recursively and come with a natural variable order (i.e., choose elements in sequential order), which facilitates the construction of the DDs. \cite{cire2013multivalued} first applied relaxed DDs to solve sequencing problems. Their implementation updates the DD during branch-and-bound search to compute more accurate dual bounds. The authors tested their DD dual bounds with satisfactory results in several sequencing problems, including the traveling salesman problem (TSP) with time windows, TSP with precedence constraints, and sequencing problems with makespan and total tardiness minimization. This work was extended by \cite{kinable2017hybrid} to tackle TSP variants with sequence-dependent cost and by \cite{castro2019mPDTSP} for pick-up and delivery problems. 

Several researchers have studied the quality of the DD bounds for sequencing problems and compared them to those obtained from a linear programming (LP) relaxation. \cite{hooker2017job} formalized some of the main components to create relaxed DDs for general sequencing problems and presents preliminary results on the bound quality for different DD construction procedures. \cite{van2018multi} used this framework to create bounds for the multi-machine scheduling problems with encouraging performance. Similarly, \cite{maschler2018multivalued} studied the DD bound quality for a prize-collecting sequencing problem and compared the bounds given by a top-down and an iterative refinement construction scheme.  Lastly, \cite{castro2018relaxed,castro2019relaxedbdds,castro2020SolvingDF} explored different DD-based relaxations for AI planning problems and show encouraging results when comparing the DD dual bounds with those from an LP relaxation.

We note that relaxed DDs and dual bound computation are strongly related to state-space relaxations \citep{christofides1981}. As we discussed in Sections \ref{sec:pre_recursive} and \ref{sec:exact_modeling}, DDs can be seen as compact representations of the state-space of a recursive model and, therefore, relaxed DDs can be interpreted as state-space relaxations. These state-space relaxations are widely used in the vehicle routing literature \citep{baldacci2011new,baldacci2012} and in other applications, e.g., cutting stock problems \citep{de2021arc}, to compute stronger bounds than LP relaxations. The main difference between these two techniques is that state-space relaxations are generally built by defining a mapping between the original states and the relaxed space. In contrast, the DD construction procedures are much more flexible since they can dynamically modify the DD to suite the user needs (see Section \ref{sec:approx_construction}).

With this in mind, one key benefit of DDs is that they can be dynamically modified to generate better bounds. The simplest alternative is to increase the width limit. However, larger relaxed DDs require more computational resources, i.e.,  memory and time for compilation. Further, empirical results show that bound improvements decrease as the DD width increases \citep{bergman2014optimization,Bergman2016book,castro2020SolvingDF}. The question of which DD width will lead to computationally efficient relaxed DDs that provide informative dual bounds is still open.

\subsubsection{Lagrangian Bounds.} \label{sec:approx_lagr}

An alternative for computing better dual bounds is to enhance the DD relaxation with dual information from an IP formulation. Consider a minimization problem with  feasible set $\feasible \subseteq \Z^n$, a linear objective function $\bc^\top\bx$, and a relaxed DD $\dd$ that over-approximates  $\feasible$, i.e., $ \feasible \subseteq \sol(\dd)$. The idea is to consider a set of $m$ valid inequalities $A\bx \leq \bb$ as penalties to the objective function to avoid paths in  $\dd$ that are infeasible. The resulting problem is a Lagrangian sub-problem
\[ \lagsubprob(\blambda) = \min\{ \bc^\top\bx + \blambda^\top(A\bx - \bb): \; \bx \in \sol(\dd)  \}, \]
where $\blambda\in \R^m_+$ are the Lagrangian penalties. 

Since $\sol(\dd)$ can be reformulated as a network flow model $\flow(\dd)$, the theoretical results of Lagrangian duality hold for $\lagsubprob(\blambda) $ \citep{ConfortiIPbook,fisher2004lagrangian}. In particular,  $\lagsubprob(\blambda)$ for any $\blambda\in \R^m_+$ is a valid dual bound for the original problem. Then, the Lagrangian dual problem seeks $\blambda$ that gives the tightest dual bounds. In our minimization example, the Lagrangian dual maximizes the sub-problem $\lagsubprob(\blambda)$ and, thus, returns a bound that is equal to or stronger than the DD relaxation:
\[ \min\{ \bc^\top\bx : \; \bx \in \sol(\dd)  \} \leq \max\{ \lagsubprob(\blambda): \; \blambda \in \R^m \}.
\]

\cite{bergman2015lagrangian} first proposed Lagrangian duality as a mechanism to enhance DD relaxations. Their approach considers a DD that represents a subset of the constraints of a problem. The Lagrangian procedure introduces dual information for the remaining constraints as penalties in the objective function. The authors tested their procedure over the TSP with encouraging results, where the DD-based Lagrangian relaxation returns significantly tighter bounds than the pure DD relaxation. \cite{hooker2019improved} further explored this idea for a family of sequencing problems and presents a detailed experimental evaluation with similar conclusions.

\cite{castro2019mPDTSP} employed this procedure for a pick-up and delivery problem and experimented with different DD and Lagrangian relaxations. Their experiments show that the DD-based Lagrangian bounds are affected by the relaxed inequalities and the constraints that are prioritized inside the DD. The authors conjectured that we can get better bounds by dualizing inequalities that are hard to represent inside the DD and prioritizing the remaining constraints in the DD relaxation.

\cite{tjandraatmadja2020incorporating} explored DD-based Lagrangian bounds for general ILP problems and proposed a decomposition approach. Instead of representing the full problem with a relaxed DD, the authors considered sub-structures that are known to be easily representable by a DD (e.g., conflict graphs). Thus, their procedure creates a DD for a specific sub-structure and uses Lagrangian duality to penalize the remaining constraints in order to improve the bounds and remove infeasible arcs. The main advantage of this procedure is that it can be used by any ILP model as long as there is a sub-structure that can be efficiently exploited by a DD.

While most works in the literature use Lagrangian duality over a single DD, 
\cite{bergman2015improved} proposed a Lagrangian decomposition approach to communicate information over multiple DDs. The idea is to have several DDs representing different constraints and use Lagrangian penalties to synchronize their solutions. The authors introduced this technique in the CP literature to improve propagation across multiple DDs. 
\cite{lange2021efficient} explored a similar idea for integer linear programming models where each linear constraint is represented with a DD. The authors employed a message-passing algorithm to solve the Lagrangian dual problem and showed that their approach obtains competitive bounds when compare to commercial solvers.

We note that despite the simplicity of this Lagrangian dual bounds, the works employing this procedure are very limited. The main advantage of Lagrangian procedures is that they avoid the construction of huge relaxed DDs to create tight bounds. As pointed out by several authors
\citep{tjandraatmadja2020incorporating, castro2019mPDTSP}, it could be more beneficial to create 
an exact or relaxed DD that represents a subset of the constraints and introduce the remaining constraints as dual penalties. By doing so, we can exploit problem structures with potentially smaller DDs and obtain better bounds than when we try to represent the complete problem with a large relaxed DD. 

\subsubsection{Primal bounds.}

In contrast to relaxed DDs, restricted DDs compute primal bounds and provide feasible solutions for a problem. \cite{bergman2014heuristic} first introduced this graphical structure as a general procedure to heuristically generate solutions for discrete optimization problems. Their approach computes primal bounds by optimizing the restricted DD using the longest/shortest path procedure for optimization problems (see Section \ref{sec:approx_dual}). Since all the paths represent feasible solutions, an optimal path corresponds to the tightest primal bound given by the restricted DD. 

The main advantage of restricted DDs is that they encode a set of feasible solutions, so they can potentially compute stronger primal bounds than other methodologies. For example, \cite{bergman2014heuristic} showed that DD primal bounds are competitive to the ones provided by IP solvers for set covering and set packing problems. Moreover, we can introduce restricted DDs into search algorithms to obtain stronger primal bounds. For instance, \cite{o2019decision} created restricted DDs to solve a TSP problem with pick-ups and deliveries in an online setting. The authors used small-width restricted DDs within a branch-and-bound search to find high-quality solutions in a few seconds.

We also note that relaxed DDs can also provide primal bounds. For example, \cite{horn2019decision} used a limited discrepancy search procedure guided by relaxed DD dual bounds to find feasible solutions for a prize-collecting job sequencing problem. Alternatively, we can extract feasible solutions from a relaxed DD using some heuristic methods that traverse the DD from root to terminal node. However, the procedure might be unsuccessful if we select a partial path that leads to infeasible solutions. 

Lastly, the primal bound of a restricted DD can also certify the feasibility of a problem. For instance, \cite{kell2013mdd} used restricted DDs to show the feasibility of multidimensional bin packing problems. If the restricted DD has at least one path, then the problem is feasible. This simple procedure proved feasibility for several instances that IP and CP technologies could not in a given time limit. 

\subsubsection{Branch-and-Bound Procedures.} 

One of the main advantages of relaxed DDs is that they can provide stronger dual bounds than an LP relaxation \citep{bergman2014optimization}. Thus, replacing the LP relaxation with a DD relaxation in a branch-and-bound procedure can be very advantageous. 

\cite{bergman2016discrete} proposed a general branch-and-bound scheme where relaxed and restricted DDs provide dual and primal bounds, respectively. The main difference with standard LP-based branch-and-bound is that their procedure branches over nodes of a relaxed DD instead of variable-value assignments. Thus, the DD-based branch-and-bound can potentially generate fewer sub-problems since each branching decision fixes the values of multiple variables at a time. \cite{gillard2021improving} presented two pruning techniques to further improve the DD-based branch-and-bound procedure that leverage local relaxed DD bound information to avoid exploring nodes that lead to sub-optimal solutions.

\cite{bergman2014parallel} extended the DD-based branch-and-bound procedure for parallel computing, where every core is responsible for a DD sub-problem. Their procedure takes advantage of the flexibility of DDs to efficiently process the sub-problems and communicate bounds between each sub-problem. 

\cite{gonzalez2020integrated} proposed a mechanism that integrates IP into the DD-based branch-and-bound. The authors modified the DD-based branch-and-bound of \cite{bergman2016discrete} so that relaxed nodes can be either solved by an IP solver or follow the original DD branching mechanism. To identify which node should be solved by an IP solver, the authors implemented a supervised learning technique that chooses nodes during search. This novel DD-based branch-and-bound procedure can be applied to any combinatorial problem and shows promising results in the maximum independent set problem and the quadratic stable set problem \citep{gonzalez2020bdd}.

DD bounds can also be integrated into other search procedures, such as a standard branch-and-bound \citep{cire2013multivalued,kinable2017hybrid,castro2020SolvingDF,tjandraatmadja2020incorporating}. The main disadvantage is that the search algorithm might not leverage the DD structure as the specialized DD-based branch-and-bound does. To take advantage of the DD structure it is important to branch on variables following the ordering in the DD. Some authors also note that a depth-first search strategy is preferable for DDs since the branching updates can be done more efficiently by removing arcs in the last branched layer \citep{cire2013multivalued,castro2019mPDTSP}. 

\subsection{CP Propagation} \label{sec:approx_propagation}

As reviewed in Section \ref{sec:exact_modeling}, several researchers in the CP community encode global constraints using exact DDs. However, the size of an exact DD grows exponentially with the number of variables, so exact DDs become impractical for large combinatorial structures. \cite{andersen2007constraint} proposed to represent global constraints with relaxed DDs to avoid the exponential size of exact DDs. 

While relaxed DDs are more flexible than exact DDs in terms of memory usage, relaxed DDs do not achieve generalized arc consistency (GAC) in polynomial time as in the case of exact DD. Since some paths in a relaxed DD are infeasible, checking if there exists a feasible solution for a specific variable-value assignment is not a trivial task. \cite{andersen2007constraint} proposed a new consistency measure to better analyze the propagation capabilities of a relaxed DD. We say that a relaxed DD $\dd=(\nodes,\arcs)$ achieves \textit{DD consistency} for a global constraint $\mathcal{C}$ if for each arc $a\in \arcs$ there exists a path $p$ that traverses $a$ and is feasible with respect to $\mathcal{C}$. Intuitively, DD consistency asks for a relaxed DD with no infeasible arcs, which is NP-hard in many cases. Note that identifying all infeasible arcs also results in identifying all infeasible variable-value assignments, thus,  DD consistency implies GAC \citep{andersen2007constraint}.

Researchers have proposed sophisticated propagation mechanisms for different global constraints to achieve DD consistency in polynomial time. A propagation procedure for a relaxed DD is defined by the set of relaxed states and the set of conditions to identify infeasible arcs (see Section \ref{sec:approx_construction}). \cite{andersen2007constraint} proposed the first relaxed DD propagators for \texttt{Linear} and \texttt{All-different} constraints. Later,  \cite{hadzic2008propagating} extended the propagator from linear inequalities to separable inequalities and showed that it achieves DD consistency in polynomial time.

Since the work of \cite{andersen2007constraint}, several papers have introduced DD propagators for other well-known global constraints. \cite{hoda2010systematic} explored DD propagators for the \texttt{Among} and \texttt{Element} constraints. The authors also proposed new conditions to identify infeasible arcs for the \texttt{All-different} propagator, improving its propagation capabilities. \cite{bergman2014mdd} presented a DD propagator for the \texttt{Sequence} constraint and show that establishing DD consistency is NP-hard. More recently, \cite{perez2017mdds}  created an DD encoding for the \texttt{Dispersion} constraint  where the DD enforces the mean value constraint and uses a cost-based propagation for the variance restriction.

Besides the encoding of existing global constraints, relaxed DDs are also a building block to create new global constraints. \cite{cire2012mdd} introduced a global constraint for disjunctive scheduling based on a relaxed DD. Their DD represents the set of job sequences and considers release times, deadlines, precedence relations, and sequence-dependent set-up times. 

Two recent works develop new probabilistic global constraints using DDs. \cite{perez2017mdds} introduced a Probability Mass Function (PMF) constraint where a DD encodes the linear inequality restricting the mean value of the variables and a cost-based propagator to ensure that the probability of every feasible assignment is inside a pre-defined range. The authors extended this constraint to consider probability distributions given by a Markov chain process \citep{perez2017mddsampling}.

DDs can also improve the propagation capabilities of a CP solver by sharing information. \cite{hadzic2009enhanced} first studied this idea using compatibility labels between multiple DDs. Their procedure constructs an exact or relaxed DD for each constraint of the problem following the same variable ordering. It then traverses the nodes of each DD to identify compatibility with nodes in different DDs. The authors showed that nodes that do not have any compatibility label can be removed since the solutions traversing that node are infeasible in all other DDs. \cite{bergman2015improved} introduced a different procedure to communicate information between DDs based on Lagrangian decomposition, which we discussed in Section \ref{sec:approx_lagr}.

While relaxed DDs can model a wide variety of global constraints, there could be alternative procedures that are more suitable for some constraints. Specifically, \cite{andersen2007constraint} noted that there are polynomial-time algorithms to enforce GAC for some constraints but it could be NP-hard for polynomial-size DDs to do so. A simple example is the \texttt{All-different} constraint that achieves GAC in polynomial time by representing the constraint as a matching problem in a bipartite graph. However, GAC is NP-hard in a relaxed DD because the GAC problem reduces to a Hamiltonian path problem.

To summarize, relaxed DDs are attractive alternatives to construct global constraints due to their propagation capabilities and flexibility to represent a wide range of combinatorial structures. However, one of the main practical limitations of relaxed DDs is their implementation, which usually requires the user to develop all the necessary components to build and propagate the DDs. Recent work by \cite{gentzel2020haddock} presented \texttt{Haddock}, a declarative language and architecture for DD compilation. The software supports a wide range of existing DD propagators and can declare multiple DDs within a CP model. Moreover, \texttt{Haddock} has comparable performance when compared to dedicated MDD propagators for different constraints.

\section{Conclusions and Future Challenges} \label{sec:conclusions}

This paper reviews recent advances using decision diagrams (DDs) to model and solve discrete optimization problems. We describe several procedures that benefit from a graphical representation of the solution set and show how these techniques can be integrated into other optimization paradigms (e.g., IP and CP). In particular, we distinguish six different ways to employ DDs: modeling, feasibility checking, solution extraction, computing primal and dual bounds, inference and propagation, and solution-space analysis. 

There are two key advantages of DDs that explain their success in many of the applications discussed. First, DDs can be efficiently optimized using different linear objective functions. This property is crucial for decomposition techniques since DDs are iteratively used, for example, in cutting plane procedures to solve the separation problem for each new fractional point. Similarly, DDs are computationally efficient when we need to extract information from the diagram multiple times, as in the case of the no-good inference procedures. 

The second key advantage of DDs is their versatility when solving optimization problems. For example, we can construct a relaxed DD to obtain a dual bound and have a heuristic procedure that extracts feasible solutions from the DD (i.e., primal bounds). Similarly, we can apply the same DD to generate valid inequalities, prune sub-optimal solutions, and infer no-good assignments for the same problem. This property makes DDs a strong and flexible optimization tool that can benefit from techniques developed in different fields. 

There are, however, several DD limitations and challenges. Most importantly, DDs can grow exponentially with the number of variables, which significantly limits their applicability. While this limitation can be partially addressed by employing approximate DDs, the quality of the approximation also depends on the size of the problem and other parametrizations. Other drawbacks include implementation challenges (i.e., lack of general and reliable DD packages for optimization) and modeling limitations (e.g., continuous variables).  

While recent works have significantly expanded the use of DDs to solve discrete optimization problems, there are many avenues yet to be explored. In particular, it is still unclear when DDs should be used. As previously mentioned, the size of a DD is a significant limitation, so researchers usually focus on applications that are combinatorially challenging but where the number of variables is small enough that the solution set can be represented with the DD. Sequencing problems are good examples since they usually have weak LP relaxations due to big-M constraints and problems with few variables can be challenging for commercial solvers \citep{cire2013multivalued,kinable2017hybrid, castro2019mPDTSP}. 

Nonetheless, we know very little about DD representability and how to efficiently construct these diagrams. We believe that it is crucial to understand which combinatorial structures are better suited to DDs so we can exploit them in our implementations. A good example of this idea is the paper by \cite{tjandraatmadja2020incorporating}. The authors proposed a DD approach that can be integrated into an IP solver, but instead of modeling the full problem with a DD they focus solely on conflict graphs, a combinatorial structure that can be efficiently model with DDs. Since representing the full problem with a DD is often impractical, knowing which sub-structures are better suited for DDs can provide a better idea of how to employ DDs in practice.  

With this idea in mind, we believe that combining DDs with existing optimization solvers can lead to state-of-the-art performance. Recent papers have shown the potential of this integration, for example, by employing DDs as a component of a decomposition algorithm (e.g., column generation, cutting planes, and Lagrangian relaxation). However, most of these works rely on existing decomposition algorithms and use DDs to efficiently solve combinatorial problems. A promising research direction is to explore new decomposition mechanisms that can leverage DDs to their full potential. DD-based branch-and-bound \citep{bergman2016discrete,gonzalez2020integrated} is a good example of a methodology specially made for DDs since the branching and bounding algorithms are designed to take advantage of the graphical structure. 

Another interesting research direction is to develop algorithms that can leverage information from multiple DDs. The main advantage of this idea is that we can represent larger problems by using multiple DDs that model small portions of them. This idea has been briefly studied in the literature, for example, with Lagrangian decomposition \citep{bergman2015improved,lange2021efficient} and network flow models with linking constraints \citep{lozano2020consistent}. However, the use of multiple DDs is still very limited, and, to the best of our knowledge, there are no efficient algorithms to integrate or share information between two or more DDs. 

Decision Diagrams are flexible optimization tools that have shown state-of-the-art results solving a wide variety of problems. DDs can be integrated with other optimization paradigms and used in different ways to solve a problem. Nonetheless, we still know very little of these graphical structures and how to properly exploit them for optimization purposes. Thus, there is a wide range of possibilities yet to be explored that can leverage DDs to solve challenging optimization problems.

\ACKNOWLEDGMENT{%
We thank the editors and reviewers whose valuable feedback helped improve the paper. Funding was provided by Agencia Nacional de Investigaci\'{o}n y Desarrollo de Chile (Becas Chile) and the Natural Sciences and Engineering Research Council of Canada (Grant Nos. RGPIN-2015-05072, RGPIN-2020-06054, RGPIN-2020-04039).
}

%
%
%


\bibliographystyle{informs2014} 
\bibliography{dd_survey_postprint} 


\end{document}